\documentclass[aip,cha,preprint]{revtex4-2}
\usepackage{amsmath,amssymb,graphicx,url}

\newcommand{\RR}{\mathbb{R}}
\newcommand{\NN}{\mathbb{N}}

\newtheorem{Theorem}{Theorem}
\newtheorem{Proposition}[Theorem]{Proposition}

\begin{document}
\title{Statistics and modelling of order patterns in univariate time series} 
\author{Christoph Bandt}
\affiliation{Institute of Mathematics, University of Greifswald, 17487 Greifswald, Germany} 
\email{bandt@uni-greifswald.de}
\date{\today}
%\subjclass[2010]{}

\begin{abstract}
Order patterns apply well to many fields, because of minimal stationarity assumptions. 
Here we fix the methodology of patterns of length 3 by introducing an orthogonal system of four pattern contrasts. These contrasts are statistically independent and turn up as eigenvectors of a covariance matrix both in the independence model and the random walk model. The most important contrast is turning rate. It can be used to evaluate sleep depth directly from EEG data. The paper discusses fluctuations of permutation entropy, statistical tests, and the need of new models for noises like EEG. We show how ordinal stationary processes can be constructed without any numerical values. An order by coin-tossing is a natural example. Every partially stationary probability measure on patterns of length $m$ can be extended to a stationary measure on patterns of infinite length. 
\end{abstract}
\maketitle

{\bf  Order patterns apply well to many fields, in particular to large data. Their importance grows as long as the size of data grows. For various questions, order patterns form a transparent alternative to machine learning.  This requires a strong theory.
Connections with dynamical systems have been thoroughly studied. 
Our paper deals with the statistics of noisy data in case of a single time series. Only patterns of length up to 5 are of interest, but the delay between time points can vary. For length 3 we define four contrasts, that is, differences of pattern frequencies. See figure \ref{ordfct}. When we require algebraic orthogonality and statistical independence of the parameters, this solution is unique.

For EEG data, a case study shows that only the turning rate contrast is relevant. As a function of time and delay, it has the power to classify sleep stages and uncover an infra-slow rhythm of the brain. 

In contrast to the values of a series, pattern frequencies in principle follow a normal distribution. This has led to asymptotic test statistics for various statistical tests. For serial dependence, permutation entropy seems a good choice because of its slim distribution under the null hypothesis.  EEG shows serial dependence at all times and scales, and so requires new models of noise. The last part of this paper sketches an approach to determine ordinal noises without using numerical values. }

\section{Introduction}\label{intro}
\subsection{A personal remark.}
Karsten Keller asked me to give the opening talk of the Orpatt22 workshop in Dresden. This was a great honour which encouraged me to tie together various lose ends of my research in recent years. Karsten and I collaborated on topology of fractals in the 90s, a work which culminated in Karsten's habilitation \cite{Ke}, and after 2000 with Bernd Pompe and others on ordinal time series analysis. I went to Dresden only for two nights, due to bad health of my wife.  We had fruitful discussions on both evenings which I shall never forget. A few weeks later, just four days after my wife, Karsten died - unexpectedly, a shock for all of us. I write this paper with great sadness and gratitude, remembering Karsten's very original, kind, constructive and optimistic way to teach and do research. 

\subsection{Contents of the paper.}
Order patterns have been applied successfully to various biomedical, geophysical, climate, finance and other data, see this issue and \cite{AKK,ZZRP}.  There are various links to the theory of dynamical systems \cite{Ami,AE,GK}, but few to the classical theory of time series, statistics and stochastic processes \cite{BS,SK,BDS}, with a recent trend to multivariate series \cite{sch,SD,NS}. Nevertheless, basic problems have to be clarified, statistical methods worked out and new models built, even in the univariate case, to establish order patterns firmly as a mathematical subject. The papers \cite{els,sou,wei} and the present note go into this direction.

\begin{figure}[h!] 
\begin{center}
\includegraphics[width=0.75\textwidth]{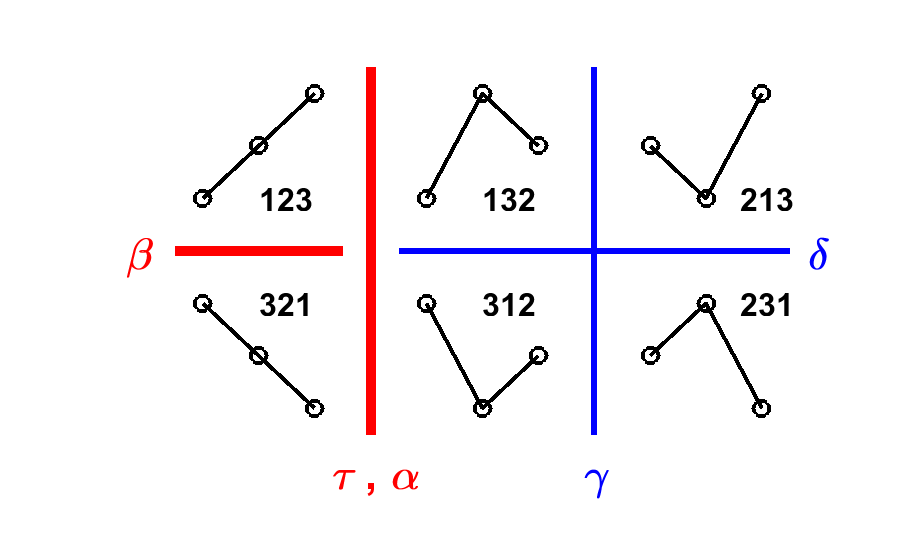}
\end{center}
\caption{The four pattern contrasts subdivide the six patterns in the best possible way. For details, see section \ref{orfunc}.}\label{ordfct}
\end{figure}  

Statistical aspects are introduced by a simulation of sample fluctuations of permutation entropy in section \ref{V3}.  Then we study order patterns of length 3. Certain pattern contrasts - that is, differences of pattern frequencies - were defined as autocorrelation functions in \cite{BS,ba14,ba15,ba17}. They are symbolically represented in figure \ref{ordfct} and discussed in section \ref{orfunc}.  In section \ref{indi}, it is proved that these differences form a canonical system of orthogonalized  patterns. They appear as eigenvectors in the random walk model as well as in the independence model. Thus they should become standard methodology for patterns of length 3. One of the pattern contrasts, the turning rate $\alpha ,$ was already studied in 1876, as noted below. 

In the second, applied part of the paper, we present a striking application of $\alpha$ to sleep stage classification by two EEG channels. Section \ref{sleepi} also discusses an apparent infra-slow biorhythm uncovered by turning rate.  We study the effects of low-pass and high-pass filtering on the order structure.  The  statistics of entropy fluctuations is used to test serial dependence of the EEG data in section \ref{serdep}.  EEG has dependencies on every scale. It requires new models.

In the third, more theoretical part, we directly construct ordinal processes without referring to numerical values. In section \ref{V6}, a natural example is constructed by coin tossing. The last section \ref{V7} defines the space of random orders on the positive integers, following a construction of Keller and Sinn \cite{KS,KSE}, and stationary probability measures on this space. An extension theorem says that any measure on permutations of length $m$ can be extended to a stationary measure on all random orders on positive integers. Various open questions are scattered throughout the paper.

\subsection{A historical note.}
The very first study of order patterns has become forgotten. In 1876, Bienaym\`e \cite{bien1} studied the `turning rate´, the relative number of local maxima and minima,  in a time series of $T$ observations.  He claimed that the turning rate is approximately normally distributed, with mean $\frac23$ and variance $\frac{8}{45 T}.$ His proof was not written down, but presented in the seminar of the Society by J. Bertrand \cite[p.225]{bien1}. However, the correct value of the variance could only be obtained by studying the 120 order patterns of length 5, in a similar way as in \cite{els,sou,wei,ba20}. Moreover, the paper contained numerous calculations with astronomical and other data which supported the normal approximation.

Today we know that the statistics of  the turning rate is exact only for white noise. We can use Bienaym\`e's discovery for a test of serial dependence, see section \ref{serdep}. The difference between independent and random (`non class\`ees arbirairement´) was not yet known, and Student's t-test was published 30 years later. This study of maxima and minima in a time series came too early.

\subsection{Order patterns and data science.}
Around 2000, the computer revolution changed society and science.  Before 2000, roughly speaking, scientists had no computers, but a lot of time to think. They developed wonderful theories. Applications were approached by simplification. Datasets were rather small. They consisted of carefully determined numbers. This time will not come back.
 
Since 2000, the growing power and availability of computers allows to get insight into very complex phenomena. Big data are everywhere, they combine numbers with structural elements, links or networks. Scientists have little time to think, and most data are not carefully checked.  Ambitious applications are taken on: 
not by theory, but by computers and automated algorithms. 
Order patterns fit well into this environment. The introduction of permutation entropy 
\cite{BP} came at the right time.

\subsection{Why patterns work so well}
Nonlinear distortion of the data, low-frequency perturbations, and outliers have little influence on order pattern methods. Intervals of missing data do not matter: already Bienaym\`e collected his turning rate from several distinct short time series. No preprocessing of data is required. Order methods are simple, fast, and easy to implement. 

Modern sensors measure in tricky indirect ways. Fluctuations of voltage indicate brain activity, carbon dioxide concentration or sea level. Changes in the environment are adjusted automatically, as a camera adapts to changing light.  As a consequence, comparison of observed values is often more meaningful than the actual values themselves. 

The most important advantage, however, is that \emph{minimal stationarity assumptions} are required.  Statistics is always based on sampling. In time series analysis, we have to sample over time points. This requires that the rules of the game do not change with time. 
Assuming a stationary Gaussian process for instance  means that all multidimensional dependencies among values reduce to correlations and do not change in time. This is never fulfilled in practice, so that perfect theoretical methods need not give nice results.  Assumptions for order patterns are less restrictive and more likely true for data.

\subsection{Work on statistics of patterns}
Having stressed practical merits of order patterns, we now regret the lack of a theory. There are some beginnings. Bandt and Shiha studied order pattern frequencies in Gaussian processes. Sinn and Keller \cite{SK} were the first to study estimators for pattern frequencies. They proved strong consistency and asymptotic normality for Gaussian processes. Their work was extended by Betken et al \cite{BDS} to long-range dependent processes and by Wei{\ss} \cite{wei} for independent variables. Tests for serial dependence, and related tests for $m$-dependence and structural breaks, were suggested by many authors, see \cite{MM,els,sou,wei,MW} and the literature quoted there. On the whole, however, classical stochastic processes and order patterns do not fit well.  Probably we must go our own way and develop a combinatorial time series analysis. Some ideas for such a theory, connected with work of Keller and Sinn \cite{KS,KSE}, will be sketched in section \ref{V7}.

\section{Pattern frequencies} \label{V2}
\subsection{Basic definitions.}
A permutation $\pi$ is an ordering of $m$ objects. In our case, the objects are consecutive time points. It is natural to start with the mathematical definition that a permutation is a one-to-one mapping from the set $\{ 1,2,...,m\}$ onto itself. A time series is also a mapping,  assigning to time points $t=1,...,T$ values $x_t=x(t).$ 
We write $\pi_k$ for  $\pi(k)$ and $\pi = \pi_1\pi_2...\pi_m.$ Thus $\pi=231$ means $\pi(1)=2, \pi(2)=3, \pi(3)=1.$ The time series 2,3,1 shall represent the permutation 231.

\begin{figure}[h!] 
\begin{center}
\includegraphics[width=0.78\textwidth]{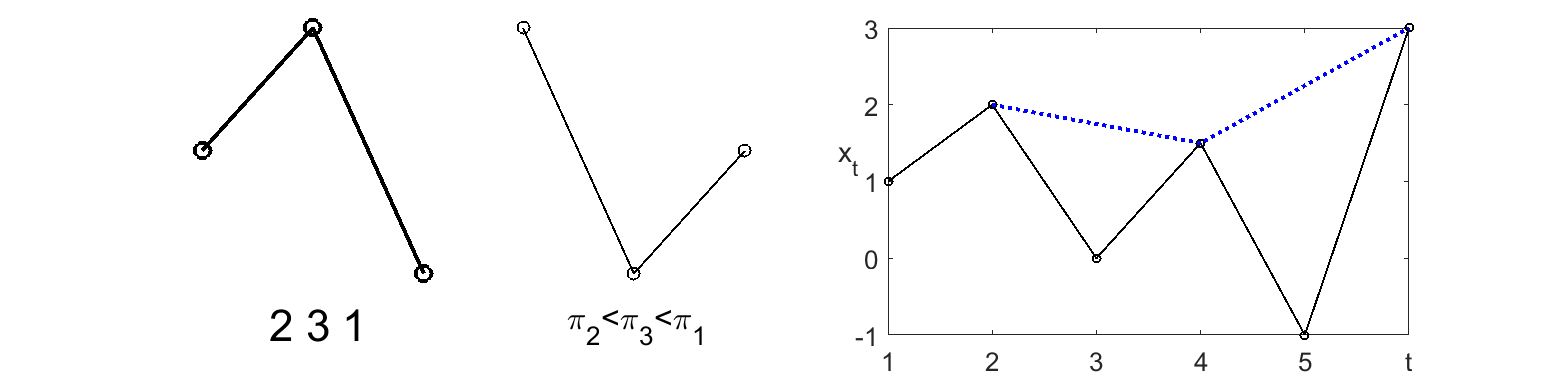}
\end{center}
\caption{Left: The graph of the permutation 231. Middle: The permutation which fulfils $\pi_2<\pi_3<\pi_1$ is called 231 by some authors. In our notation, it is the inverse 312 of 231. Right: Example time series. The dotted line indicates a pattern with $d=2.$}\label{fig1}
\end{figure}

The graph of a permutation $\pi$ represents an order pattern, see figure \ref{fig1}. Any pattern that shows the same comparisons between a series of $m$ numbers or objects is considered equivalent to $\pi .$ Thus $(0, 6,  -2)$ and $(0.199, 0.2, -4)$ are equivalent to 231. We consider $\pi$ as the name of an equivalence class of order patterns, and as standard form, in a square with uniformly spaced numbers.

Most authors write 231 for the permutation which fulfils $\pi_2<\pi_3<\pi_1.$  They read the values on the $x$-axis from bottom to top while we read them on the $y$-axis from left to right. Note that this notation is counterintuitive since the time series $(3,1,2)$ then represents the pattern 231, and the time reversed series $(2,1,3)$ does not correspond to the time-reversed pattern 132. However, figure \ref{fig1} shows the only case for $m=3$ where the two notations differ.  Some authors also read from top to bottom \cite{BDS}.

With $x=x_1,x_2,...,x_T$ we denote a time series of length $T.$ 
For $n$ time points $t_1<t_2<...<t_m$ between 1 and $T$ we consider the pattern of corresponding values. We define
\begin{equation} (x_{t_1},x_{t_2},...,x_{t_m}) \mbox{ shows pattern $\pi$ when }
   x_{t_i}< x_{t_j} \mbox{ if and only if } \pi_i< \pi_j \ .
\label{pati}\end{equation}
Actually, the pattern $\pi$ can be determined directly from the time series by calculating ranks: $\pi_k$ is the number of time points $t_j$ for which $x_j\le x_k.$ 

In this paper, we study only equally spaced time points $t< t+d< t+2d<...< t+(m-1)d.$ We call  $m$ the length and $d$ the delay (or lag) of the pattern. We always take $m<7$ and mainly consider $m=3$ and $m=4.$ However, we vary $d$ as much as possible.  

{\bf Example. } The time series  $x= 1, 2, 0, 1.5, -1, 3$ in figure \ref{fig1} shows patterns 231, 312, 231, 213 for delay 1, and 321, 213 (dotted line) for delay 2. 

Now we fix a length $m.$ In examples, we usually take $m=3.$ We determine relative frequencies of all patterns $\pi$ of length $m,$ for various delays $d.$   We divide the number of occurences of $\pi$ by the number of time points where $\pi$ could occur.
\begin{equation}
p_\pi (d) = \frac{1}{T-(m-1)d}\ \# \{ t\  |\  1\le t\le T-(m-1)d, \ (x_t, x_{t+d},...,x_{t+(m-1)d}) \mbox{ shows pattern } \pi \} \ .
\label{esti}\end{equation}
For the example above, $p_{231}(1)=1/2,\ p_{231}(2)=0$  and  $p_{213}(1)=1/4,\ p_{213}(2)=1/2.$ In a language like Matlab or R, five lines of code suffice to determine all $p_\pi(d)$ for a given $d$ \cite{ba14}.

\subsection{Assumptions and estimates.}
We are not interested in the time series for its own sake. We rather want to know the laws of an underlying mechanism which produced this time series and which could also produce many similar time series by random variation. It is custom to represent this mechanism by a sequence $X=X_1,X_2,...$ of random variables on a probability space $(\Omega , P).$
We call $X$ a stochastic process and imagine the time series $x$ is a random choice from the ensemble $X.$

Since we want to find the laws of $X$ from one single time series, we must assume that the values and dependencies of the $X_k$ do not change in time. Usually $X$ is required to be stationary. That is, the $m$-dimensional distribution of $(X_t,X_{t+1},...,X_{t+m})$ does not depend on $t,$ for all $m.$ This is a fine starting point for theory but never fulfilled in reality.
Weak stationarity says that mean $M(X_t)$ and $Cov (X_t,X_{t+k})$ do not depend on $t$ \cite{GS}. Even this need not be true in practice, for instance in EEG measurements. We make a less stringent assumption.

{\bf Order stationarity. } For a given pattern $\pi$ with length $m$ and delay $d,$ let 
\begin{equation}\label{ost}
P_\pi(d)= P\{ (X_t,X_{t+d},..., X_{t+(m-1)d}) \mbox{ shows pattern }\pi \}
\end{equation}
denote the probability that the process $X$ shows pattern $\pi$ just after time $t.$  We say that $X$ is  \emph{order stationary for patterns of length $m$} if $P_\pi(d)$ is the same value for all $t.$
This should hold for all patterns $\pi$ of length $m$ and all delays $d=1,...,d_{\rm max}.$
We take small $m,$ usually 3 or 4, and choose $d_{\rm max}$ much smaller than the length $T$ of our time series. Then there is hope that such order stationarity is fulfilled in reality, at least approximately. Note that stationary processes and processes with stationary increments (cf. next section) are order stationary for all $\pi, m, d.$ 
The second requirement assumed throughout this paper is  
\begin{equation}\label{equa}
P\{ X_s= X_t \}=0 \quad\mbox{  for }\  s\not= t\ .  
\end{equation} 
We exclude equality since order patterns become complicated in case of equal values. When they appear in reality, we can just ignore them like missing values, and adapt the denominator in \eqref{esti}. Condition \eqref{equa} can also be guaranteed in practice by adding a tiny random noise to the data.

The $P_\pi(d)$ are the parameters of the process which we would like to determine. The numbers $p_\pi(d)$ are the concrete frequencies obtained from our data series $x.$ It is natural to take $p_\pi(d)$ as estimate of the unknown probability $P_\pi(d).$ This naive estimator is unbiased \cite{SK}. For Gaussian \cite{SK,BDS} and $m$-dependent processes  \cite{sou,wei} the estimator is strongly consistent and asymptotically normal.  The estimation error, that is, the variance of the estimator,  can only be determined for more specific models.

\subsection{Independence model and random walk model}
There are essentially only two random processes for which we know the probabilities  $P_\pi(d)$ rigorously. The first one is a sequence of independent and identically distributed  random  variables $X_1,X_2,...$ which is called an i.i.d.\ process, and white noise in the case of normal distribution. This process induces the equidistribution on patterns: $P_\pi(d)= 1/m!$ for all patterns $\pi$ of length $m$ and for all delays $d.$ 
The reason is that an i.i.d.\  process is exchangeable: $(X_1,X_2,...,X_m)$ has the same distribution as $(X_{\pi_1}, X_{\pi_2},...,X_{\pi_m})$ for any permutation $\pi .$  

{\bf Note. } The one-dimensional distribution of the values, which has nothing to do with order patterns, need not be specified. For simulations we take standard normal distribution.  The class of exchangeable processes theoretically also contains mixtures of several i.i.d.\ processes. However, for an ergodic stationary process, Bandt and Shiha \cite{BS} proved that the one-dimensional distribution plus all pattern frequencies completely describe all finite-dimensional distributions of the process. Thus order structure and one-dimensional distribution are two complementary parts which together characterize a stationary process.

An i.i.d.\  process is self-similar in the sense that $(X_1,X_2,...,X_m)$ and $(X_d,X_{2d},...,X_{md})$ have the same distribution. Thus $P_\pi(d)$ does not depend on $d.$ Actually, this seems the only stationary process for which $P_\pi(d)$ is known for arbitrary large $d.$

This process is the null hypothesis for testing independence of the variables $X_k,$ see section \ref{serdep}. 
It should be mentioned that although the $X_k$ are independent, the events $X_k<X_{k+1}$ are dependent. Under the condition $X_1<X_2<...<X_k$ we have $P\{ X_k<X_{k+1} \}=1/(k+1)$ while in the mean this probability is 1/2.\medskip

The second process is \emph{symmetric random walk} where we set $X_0=0$ and let the increments $I_k=X_{k+1}-X_k$ be independent  random variables with identical distribution for $k=0,1,...$ with density function fulfilling $\varphi(x)=\varphi(-x).$ This is a non-stationary process with stationary and independent increments $I_k.$ In the case of Gaussian distribution the process is self-similar, so that $P_\pi(d)$ does not depend on $d.$ This case will be called Brownian motion.

From an ordinal viewpoint, it is not appropriate to consider random walk just as a cumulative sum of an i.i.d.\  process. For example,  geometric Brownian motion $Y_k=\exp(X_k)$ has exactly the same order structure and is not a cumulative sum. It is a standard model for data from financial markets. Gaussian random walk also represents a null hypothesis \cite{ba19}.

To determine the $P_\pi(d)$ for $d=1$ and patterns of length 3, we need only assume a symmetric random walk \cite{sou}. Then 
\begin{eqnarray} P_{123}=P\{ X_1<X_2<X_3\} =P\{ I_1>0, I_2>0\}=(\frac{1}{2})^2=\frac{1}{4}=P_{321}\ ,\notag\\
 P_{132}=P\{ X_1<X_3<X_2\} =P\{ I_1>0, I_2<0 \mbox{ and } |I_2|<I_1\}=(\frac{1}{2})^3=\frac{1}{8} \ , \label{brown3}
\end{eqnarray}
and 1/8  also for the remaining probabilities. Here we used the Markov property of random walks: the occurence of an order pattern among the $X_s$ with $s\ge t$ is independent of order patterns of the $X_s$ with $s\le t.$

For length 4, pattern probabilities depend on the distribution of the increments. For Gaussian random walk they were calculated in \cite{BS}. For length $\ge 5,$ they can be expressed as multidimensional integrals which have no analytic solution.  Elizalde and Martinez \cite{EM} studied classes of permutations which have the same pattern frequencies for random walks with arbitrary distribution.

Pattern probabilities have also been determined for moving average and autoregressive processes, mostly by simulation.  See  \cite{BS, sou} and table \ref{arpat} below.  \smallskip

{\bf Open problems. } Both stationary processes and processes with stationary increments are order stationary.  Do pattern probabilities indicate whether a process is stationary or has stationary increments? Find examples of order stationary processes which are neither stationary nor have stationary increments! (A candidate is given in section \ref{V6}.) Does order structure plus one-dimensional distribution of increments completely describe a process with stationary increments? \smallskip

{\bf Example. }  The condition $p_{12}=1$ characterizes an increasing series. It is a process with stationary increments. The one-dimensional distribution of increments must be concentrated on positive numbers, otherwise it is arbitrary.
The conditions $p_{231}=p_{312}=p_{321}=0$ describe a similar but slightly more complicated case.

\section{Permutation entropy} \label{V3}

\subsection{Permutation entropy at its maximum}\label{chisi}
Let $S_m$ denote the set of all permutations of length $m$ which usually will be between 3 and 6. The permutation entropy is defined as 
\begin{equation}\label{pe}
H=H_m(p) = - \sum_{\pi\in S_m} p_\pi  \log  p_\pi  \ . 
\end{equation}
We use natural logarithm, so $H=1$ corresponds to $1/\log 2=1.44$ bit.
The maximum value $H^{\rm max}=\log m!$ is assumed for an i.i.d.\  process and the minimum value 0 for monotone time series \cite{BP}.  The scale for $H$ is nonuniform. For a period 2 sequence, like $x_t=(-1)^t/t=-1,\frac12,-\frac13,...$ with $P_{132}=P_{312}=\frac12 ,$ we have $H_3=\log 2=0.69,$ and for period 3 with $P_{123}=P_{231}=P_{312}=\frac13 ,$ we get $H_3=\log 3=1.10.$ Thus $H$ works very well in detecting periodicities and distinguishing small periods and other dynamical phenomena  \cite{Ami}. On the other hand, noisy time series, coming from sources like EEG or geophysical measurements, all have entropy very near to the maximum. For Brownian motion, \eqref{brown3} yield 
$H_3=\frac{5}{2}\log 2 = 1.7345$ while  $H_3^{\rm max}=\log 6=1.7918.$ These are two fundamentally different processes.

Our question here is whether permutation entropy can distinguish statistically between such noisy processes. To this end, we consider white noise as our point zero. The essential quantity is not $H,$ but its deviation $H^{\rm max}-H$ from its maximum.
We define the distance to white noise \cite{ba17}
\begin{equation}
\Delta^2 =  \sum_{\pi\in S_m} \left( p_\pi  -\frac{1}{m!}\right)^2    \ . \label{d2}
\end{equation}
For simplicity, we do not include a delay parameter $d.$ We also drop the subscript $m$ since we now consider $H(p)$ as a differentiable function of $m!$ variables $p=(p_\pi)_{\pi\in S_m}$ near the white noise parameter  $p^o=(1/m!)_{\pi\in S_m}.$  Our purpose is to get an idea of the statistical fluctuations of $H$ in different samples of an i.i.d.\  process.  We need the second order Taylor expansion of $H$ at $p^o.$

\begin{equation}
H=H_m(p) \approx \log m! - \frac{m!}{2}\Delta^2   \quad \mbox{ or }\quad  \frac{m!}{2}\Delta^2  \approx H^{\rm max}-H \ .   \label{tayli}
\end{equation}
This is easy to check: $\frac{\partial H}{\partial p_\pi}=-1-\log p_\pi$ is a constant for the i.i.d.\  process. So the first-order term vanishes, as it should for a maximum. The non-zero second derivatives are 
$\frac{\partial^2 H}{\partial p_\pi^2}=-\frac{1}{p_\pi}$ with value $m!$ at $p^o.$ So the second order term is $\frac{m!}{2}$ times the sum of the $(p_\pi-1/m!)^2.$

For a check of \eqref{tayli}, take symmetric random walk with 
$H_3(p)=1.7345.$  The Taylor approximation, calculated as $H_3(p)\approx \log 6 -3 \Delta^2=1.7918- 1/16 =1.7293,$ coincides with two decimals accuracy.

\subsection{Simulated fluctuations of entropy estimates}\label{flucti}
We have proved that near the i.i.d.\  process, $H^{\rm max}-H$ is the quadratic quantity
$\frac{m!}{2}\Delta^2 .$ For all other processes, the Taylor approximation is first order and the deviations are linear.  Now we shall discuss the consequences for statistical fluctuations of $H$ in sample series of different lengths $T.$ The simulations for white noise and Brownian motion show different behavior.

\begin{figure}[h!] 
\begin{center}
\includegraphics[width=0.48\textwidth]{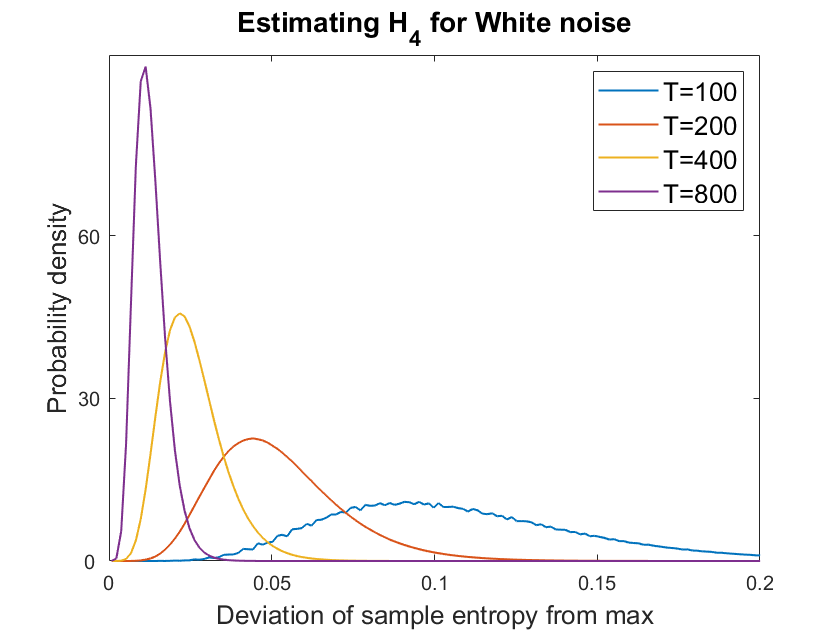}\quad
\includegraphics[width=0.48\textwidth]{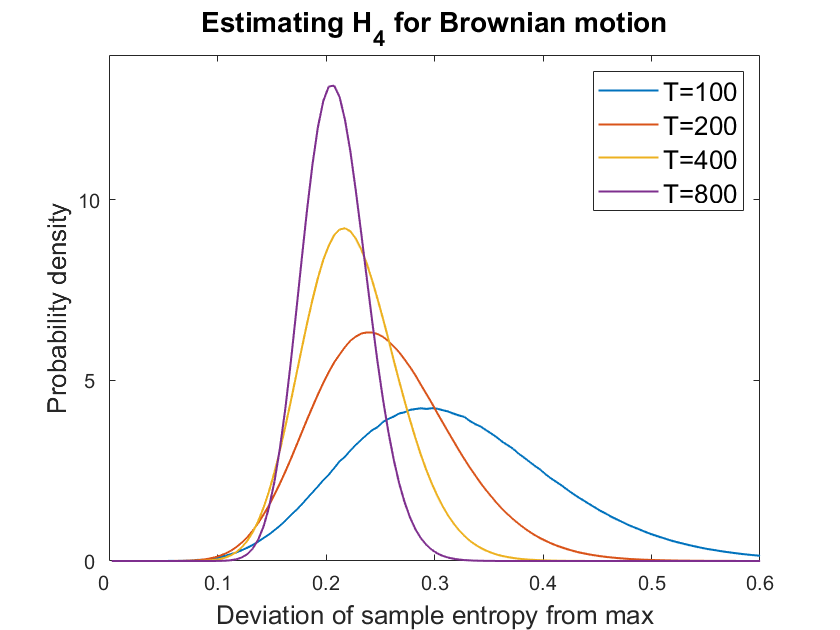}
\end{center}
\caption{Left: deviations $H_4^{\rm max}-H_4$  of 100000 white noise sample series of different length $T.$ The deviations scale like $1/T.$  Right: the same simulation for Brownian motion. The bias scales like $1/T,$ but the standard error scales like $1/\sqrt{T}.$}\label{simul1}
\end{figure}

For figure \ref{simul1}, we simulated 100000 time series of different lengths $T$ and determined the distribution of $H_4^{\rm max}-H_4.$ The true value is 0 for white noise and 0.194 for Brownian motion. This does not agree with the mean of the distributions. There is considerable bias in the statistical estimates. The bias is always to the right, to larger deviations from $H^{\rm max}.$ It scales like $1/T$ with the sample size.

For white noise the bias is easy to explain. The sample frequencies $p_\pi$ are never all equal, so the entropy of the sample series cannot reach the value $H^{\rm max}.$ Surprisingly, we have the same bias also for other processes. Table \ref{tabi} shows that the bias is divided by two when we double the sample length $T.$ This observation helps to extrapolate the true value and correct the bias.

\begin{table}[h!]
\begin{tabular}{|c|c|c|c|c|c|}\hline
$T$ &100&200&400&800&limit\\ \hline
mean&0.320&0.254&0.224&0.209&0.194\\ \hline
bias&0.126&0.060&0.030&0.015& 0 \\  \hline
std.dev. $\sigma_T$&.0974&.0648&.0437&.0304& 0 \\ \hline
$\sigma_T/\sigma_{2T}$&1.50&1.48&1.44& - &$\sqrt{2}$\\   \hline
\end{tabular} 
\caption{Mean, bias and standard deviation of $H_4^{\rm max}-H_4$ taken over 100000 sample series of Brownian motion with different $T.$ Scaling of bias by 2 is obvious, scaling of $\sigma$ by $\sqrt{2}$ follows from the central limit theorem.} \label{tabi}
\end{table}

What about the variation of the entropies of samples?  It turns out that the estimates for white noise are much more accurate than for any other process, due to its special Taylor expansion. For $T=800,$ for instance, figure \ref{simul1} shows estimates between 0 and 0.04 for white noise and between 0.1 and 0.32 for Brownian motion. Actually, the standard deviation of the distribution scales with $1/T$ for white noise and with $1/\sqrt{T}$ for other processes.

\begin{figure}[h!] 
\begin{center}
\includegraphics[width=0.48\textwidth]{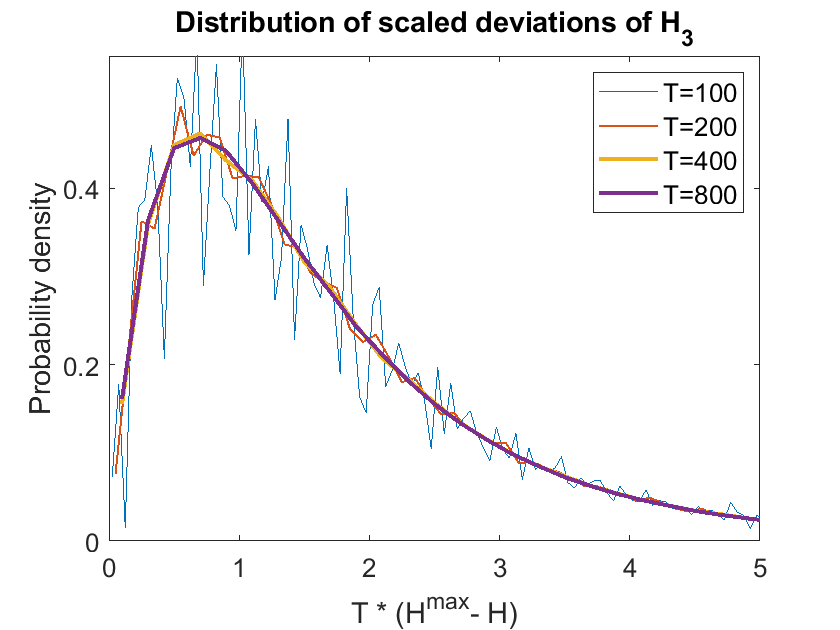}\quad
\includegraphics[width=0.48\textwidth]{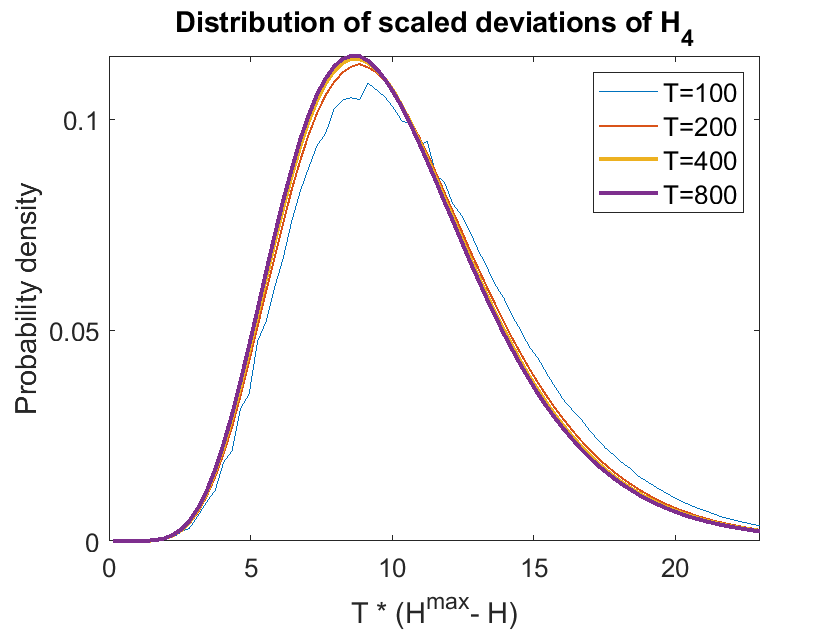}
\end{center}
\caption{Distribution of scaled deviations $Z=T*(H^{\rm max}-H)$ of 100000 sample series of different length $T$ from the white noise process verifies the scaling. The distributions do not depend much on $T$ and approach a limit function \cite{wei}. The distributions for $m=3$ show irregularities due to discrete effects.}\label{simul2}
\end{figure}  

The scaled deviations $Z=T*(H^{\rm max}-H)$ in figure \ref{simul2}
illustrate the $1/T$-scaling of the statistical fluctuations of white noise.  There is a good agreement for different $T,$ and there is a limit density for $T\to\infty .$ It is calculated in Corollary 3.2 in Weiss \cite{wei} which extends to $m>3$ when corresponding eigenvalues are determined.

For $m=3,$ and in case $m=4$ for $T=100,$ the density functions look irregular. This is not statistical inaccuracy. The density did not change when the number of simulations was increased or decreased. The density did change, however, when $T=100$ was replaced by 99 or 101. The irregularities are caused by the fact that absolute frequencies of patterns are relatively small integers - so there are not so many possible frequency combinations. To avoid such effects, just work with $m=4$ and $T\ge 200.$

\subsection{Tail probabilities and quantiles}\label{tailprob}
When we know the distribution of $Z$ in an i.i.d.\  process, in particular the $p$-values or tail probabilities $p(z)=P(Z\ge z),$ we can perform a statistical test of serial dependence. Such tests are very important for time series analysis. Whenever a model for a time series fits the data, the residues - theoretical minus observed values - should be indistinguishable from i.i.d.\  random numbers. An overview of the multitude of serial dependence tests  can be found in Matilla-Garcia and Ruiz Marin \cite{MM}.  

The paper \cite{MM} suggested to take permutation entropy, more precisely the values $2Z,$ as criterion for a serial dependence test. Unfortunately, the authors wrongly claimed that under the null hypothesis of an i.i.d.\  process, $2Z$ follows a $\chi^2$ square distribution  with $m!-1$ degrees of freedom. Actually $Tm!\Delta^2 $ fulfils a well-known formula for the $\chi^2$-test:
\[  2Z\approx Tm!\Delta^2  = Tm!\sum_\pi (p_\pi-1/m!)^2 =\sum_\pi \frac{(Tp_\pi -T/m!)^2}{T/m!} =\sum\frac{(B-E)^2}{E}\]
where $E=T/m!$ is the expected and $B=Tp_\pi$ the observed absolute frequency of pattern $\pi$ in the sample series. This roughly shows the type of limit distribution, at least for $H_3,$ when we take four degrees of freedom.  However, Elsinger \cite{els} pointed out that besides degrees of freedom, other strong corrections are necessary due to correlation and different variances of pattern frequencies, see section \ref{orfunc}. 
Following Elsinger \cite{els}, two different tests were worked out by Sousa and Hlinka \cite{sou} and Weiss \cite{wei} in this special issue. The main tool in both cases is principal axes analysis, see section \ref{indi}. Weiss \cite{wei} determines the limit distribution as a convolution of multiples of $\chi^2_1$ for $m=3.$ His Corollary 3.2 extends to $m>3$ when corresponding eigenvalues are determined. Elsinger \cite{els,sou} uses a transformation to standard multivariate normal distribution and a classical $\chi^2$ statistics. 

\begin{figure}[h!] 
\begin{center}
\includegraphics[width=0.48\textwidth]{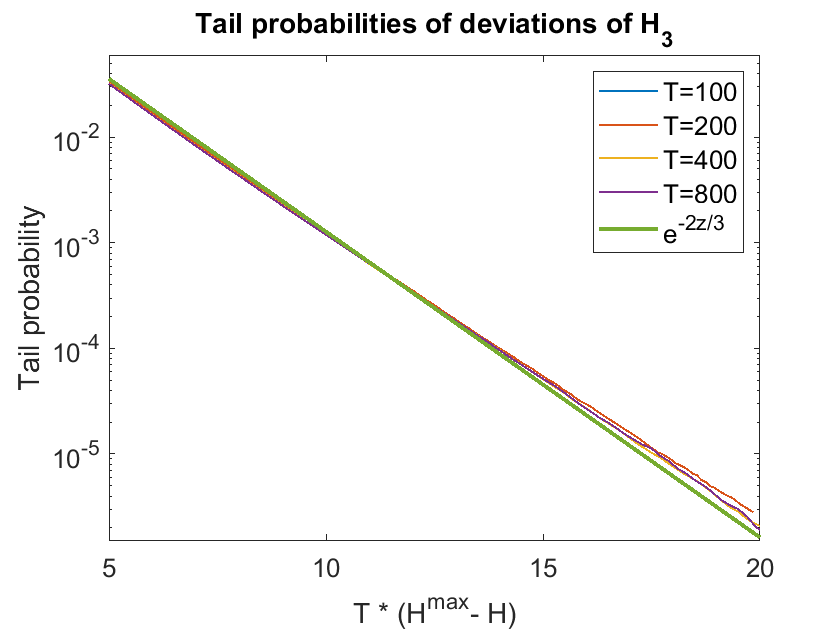}\quad
\includegraphics[width=0.48\textwidth]{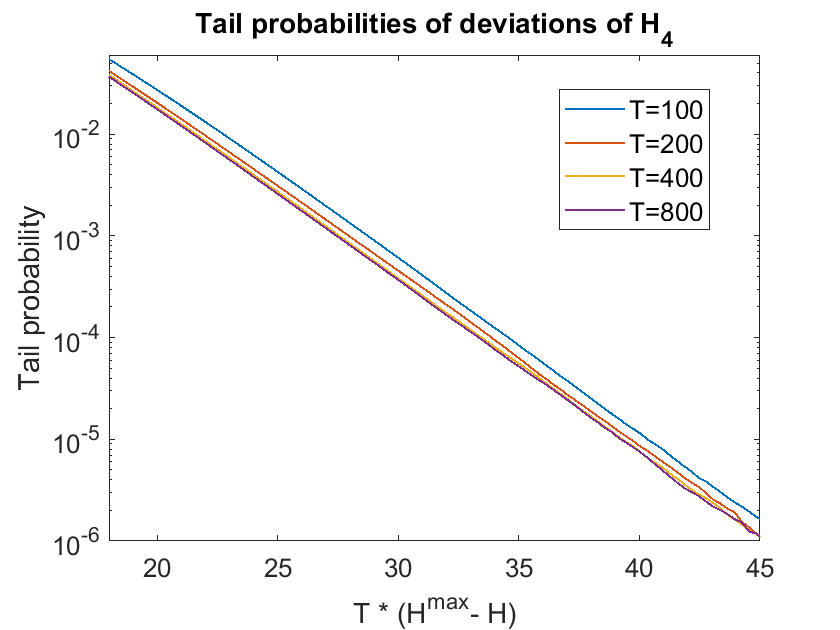}
\end{center}
\caption{To perform tests of serial dependence, we need tail probabilities ($p$-values). For $H_3,$ they almost coincide for different $T.$  The empirical formula $p=\exp(-2z/3)$ yields good $p$-values even for larger $z$ were $\chi^2$ is misleading. For $H_4,$ dependence on $T$ is small but visible up to $T=800.$ }\label{tailp}
\end{figure}

Such asymptotic results are important for our understanding. However, they do not tell us what happens for a particular $T.$  For tests of serial dependence in section \ref{serdep} we include figure \ref{tailp} and the critical quantiles ($z$-values with $P(Z\le z)=0.95,$ say) in table \ref{tailq}, obtained by direct simulation (100000 samples for each $T$). For $H=3,$ the formula $p=e^{-2z/3}$ works well even for rather large $z,$ as figure \ref{tailp} shows. The table for $H_3$ agrees with the table for $2Z$ in \cite{els}, p.21.
 
\begin{table}[h!]
\begin{tabular}{|c|c|c|c|}\hline
sample size&95\%&99\%&99.9\%\\ \hline
$T=800$&4.46&6.81&10.36\\
$T=400$&4.47&6.82&10.37\\ 
$T=200$&4.49&6.91&10.39\\ \hline
$e^{-2z/3}$&4.49&6.91&10.36\\  
$\chi^2_4(2z)$&4.74&6.64&9.24\\  \hline
\end{tabular} \qquad
\begin{tabular}{|c|c|c|c|}\hline
sample size&95\%&99\%&99.9\%\\ \hline
$T=\infty$&17.22&21.55&27.48\\
$T=800$&17.29&21.63&27.56\\
$T=400$&17.39&21.75&27.69\\ 
$T=200$&17.64&22.06&28.08\\ 
$T=100$&18.40&22.89&28.88\\  \hline
\end{tabular}
\caption{Critical quantiles of $Z$ for tests of serial dependence with  $Z=T*(H^{\rm max}-H)$ for $m=3$ (left) and $m=4$ (right). Values for $T=\infty$ were extrapolated.} \label{tailq}
\end{table}

\section{Statistics of patterns of length 3}\label{disti}
\subsection{Why short patterns?} 
There are three reasons for considering only patterns of length up to 5.
\begin{itemize} \item Estimation: statistical estimates of pattern frequencies won't be accurate for patterns of length $\ge 6.$ Already for the 120 patterns of length 5, one needs time series with several thousand values for a reasonable estimate. 
\item Interpretation: it is difficult to give a meaning to longer patterns.
\item  Assumptions of stationarity are more restrictive for longer patterns. We have no experience and intuition for such multidimensional dependencies in random processes.
\end{itemize}
In this section, we study only patterns with length $m=3$ and delay $d=1,$ so we omit these parameters. 
We take an order stationary random process and imagine we sample a lot of time series of length $T$ from this process, as in the simulation of $H.$  Our time series should never be shorter than $T=200,$ to avoid discrete effects as in figure \ref{simul2}.
.
\subsection{Properties of estimators of pattern frequencies}
It is important to understand that in contrast to the values of a time series, pattern frequencies are `always' asymptotically normally distributed. We only have to assume that the dependence between occurence of patterns at two time points $s,t$ decreases fast when the time difference $|s-t|$ increases. A special case treated in \cite{els,sou} is $m$-dependence which says that the patterns become independent if $|s-t|>m.$ This assumption is true for Brownian motion, MA-processes and for the coin-tossing model in section \ref{V6}, and it can be expected to hold in practice. Under such an assumption even the combination of all six pattern frequencies will approximately follow a multivariate normal distribution. This follows from central limit theorems, cf. \cite{sou,wei}. We shall not rely on asymptotic statements which contain no quantitative estimates. But they are good for orientation. Since correlations describe all dependencies between variables of a multivariate normal distribution, we have to study covariances of pattern estimators.

Indeed the estimators $\hat{p}_\pi,$ that is, the sample frequencies taken as functions on the space of all sample series, are not binomially distributed. They are heavily correlated, and the variance is different for different $\pi .$ For example, an increasing time series contains $T-2$ patterns 123, but the pattern 231 can appear at most $\frac12 (T-1)$ times since two occurences cannot be consecutive.  Combinatorial calculations performed by Elsinger \cite{els}, see also \cite{sou,wei} show that asymptotic variances and covariances scale with $1/T,$ and are different for different stochastic processes. The exact form of the covariances is $a/T +b/T^2,$ but the second term will not be discussed. It is negligible if $T\ge 200.$
Here are the asymptotic covariance matrices for an i.i.d.\ process \cite{els} and symmetric random walk \cite{sou} which we need later. The six rows and columns correspond to the six permutations 123, 132, 213, 231, 312 and 321. 

\begin{equation}\label{sigW}
\Sigma_W=\frac{1}{360T}\left(\begin{array}{cccccc}
 46&-23&-23&7&7&-14\\
-23&28&10&-2&-20&7\\
-23&10&28&-20&-2&7\\
7&-2&-20&28&10&-23\\
7&-20&-2&10&28&-23\\
-14&7&7&-23&-23&46\\
\end{array}\right) \end{equation}
\begin{equation}\label{sigB}
\Sigma_B=\frac{1}{192T}\left(\begin{array}{cccccc}
60&-6&-6&-6&-6&-36\\
-6&15&7&-9&-1&-6\\
-6&7&15&-1&-9&-6\\
-6&-9&-1&15&7&-6\\
-6&-1&-9&7&15&-6\\
-36&-6&-6&-6&-6&60\\
\end{array}\right)
\end{equation}

It is not difficult to see that the six frequencies of patterns of length 3 have only \emph{four degrees of freedom.} There are two constraints. The first is obvious. 
\begin{equation}  
  \sum_\pi  p_\pi =1 \ . \label{cons1}
\end{equation}
Since in any time series the numbers  of local minima and local maxima can differ by at most 1, it is true that $|p_{213}+p_{312}-p_{231}-p_{132}|\le 1/T$ for the frequencies in any time series of length $T.$ For stochastic processes, even nonstationary ones, we then have the asymptotic equation
\begin{equation}  
    p_{213}+p_{312}-p_{231}-p_{132} =0 \ . \label{cons2}
\end{equation}
For stationary processes this holds immediately, without letting $T\to\infty .$ One reason is that $P\{ X_1<X_2\}= p_{12*}=p_{123}+p_{132}+p_{231}$   and   $P\{ X_2<X_3\}= p_{*12}=p_{312}+p_{213}+p_{123}$ must be equal.

\subsection{The pattern contrasts for length 3}\label{orfunc} 
It is rare that particular patterns have a meaning, as 1324 and 1342 in the heart rate study of Parlitz et al. \cite{par}.  Usually certain differences of pattern frequencies are more powerful than the patterns themselves. They were introduced in \cite{BS,ba14,ba15,ba19} as a kind of autocorrelation functions of $d=1,2,...$ which describe the behavior of data, for instance periodicities, short-time dependence or self-similarity.   
This practically important aspect is not studied here.  We have $d=1$ and focus on showing that the pattern contrasts form an ideal orthogonalized system of patterns, both from an algebraic and a statistical viewpoint. Let us now introduce the pattern contrasts which were illustrated in figure \ref{ordfct}. `Contrast' is a basic concept of the analysis of variance, contained in every advanced textbook of statistics and in the special introduction \cite{RR}. It means `weighted sum of frequencies' and is used here in exactly the same spirit. It is actually a difference since the sum of weights is usually zero.

\begin{equation}
\mbox{The up-down balance }\qquad \beta = p_{123}-p_{321} =  p_{12}-p_{21} = 2  p_{12}-1 \label{beta}
\end{equation}
distinguishes upward and downward patterns. It is zero for Gaussian and reversible processes and positive for time series which increase slowly and decrease fast. It is a symmetry parameter, with extreme values $-1$ for a decreasing and $+1$ for an increasing time series. As \eqref{beta} shows, this difference comes from patterns of length 2. All identities among pattern frequencies hold for order stationary processes. For time series there is often a small boundary error \cite{BS,ba14}.
\begin{equation}\label{tau}
\mbox{The persistence }\qquad\qquad \tau = p_{123}+p_{321}-\frac{1}{3} =\frac23 -\alpha
\end{equation}
has properties very similar to classical and rank autocorrelation, when considered as a function of $d.$  It can also be considered as a smoothness parameter since it distinguishes straight and broken patterns. It is just the negative of the turning rate $\alpha$ of Bienaym{\`e}, plus $\frac23.$ The constant of $\tau$ was chosen so that it is zero for white noise, like the other pattern contrasts. Smoother time series have smaller turning rate and larger persistence. The maximum value $\frac23$ of $\tau$ occurs for a monotone series. The minimum value is $-\frac13$ and occurs for an alternating sequence, as mentioned after \eqref{pe}. Negative persistence is rare for $d=1,$ but typically for $d>1$ when $2d$  is near to a period of the process. $\tau$ or $\alpha$ is the most important pattern contrast, and both versions are useful. Turning rate has the advantage of clear interpretation.
\begin{eqnarray}
\mbox{The rotational asymmetry }\qquad \gamma &=& p_{213}+p_{231}- p_{132}-p_{312}\\ &=& 2(p_{213}-p_{312})=2(p_{231}-p_{132}) \notag
\end{eqnarray}\label{gamma}
distinguishes the patterns with middle rank 2 at the beginning from those where 2 is at the end. It is positive for an oscillation with increasing amplitude, as in figure \ref{fig1}, negative for damped oscillation, and zero for Gaussian and reversible processes. The extreme values $\pm 1,$ occur for alternating series.  'Rotational' refers to a geometric half-turn of the pattern. Obviously $\gamma$ measures the reversibility of the series.
\begin{eqnarray}
 \mbox{The up-down scaling }\qquad \delta &=& p_{132}+p_{213}-p_{231}-p_{312}\\
 &=& 2(p_{132}-p_{312})=2(p_{213}-p_{231}) \notag
\end{eqnarray}\label{delta}
got its name from the fact that $\delta(d)=\beta(2d)-\beta(d) .$ It is tightly connected with $\beta .$ It is zero for Gaussian, reversible, and self-similar processes. Thus the pattern contrasts $\beta , \gamma , \delta$ all  indicate deviation from symmetry and reversibility, but in different ways.  

\begin{table}[h!]
\begin{tabular}{|c|c|c|c|c|}\hline
noise&$\beta$&$\tau$&$\gamma$&$\delta$\\ \hline
normal&0&.086&0&0\\
uniform&0&.083&.042&0\\ 
Bernoulli&0&1/6&1/6&0\\ \hline
triangular&-.074& .088& .026&.048\\ 
exponential&.161&.107&.002&-.095\\ \hline
\end{tabular}
\caption{pattern contrasts for the AR(1) model $X_t=\frac12 X_{t-1}+\varepsilon_t$ with different types of noise $\varepsilon_t$}\label{arpat}
\end{table}

To give an idea of the size of the parameters, we simulated AR(1) processes $X_t=\frac12 X_{t-1}+\varepsilon_t$ with different distributions of the i.i.d.\ noise $\varepsilon_t.$
For symmetric distributions, we have $\beta=\delta=0.$ However, $\gamma$ is positive even for uniform noise although such a minor effect is not visible in small sample series. (To distinguish it from normal noise with error probability $5\%  ,$ let $T=3000 .$ Then 98\%  of the samples with uniform noise will be recognized correctly.) Bernoulli noise means $\varepsilon_t=\pm 1$ with probability $\frac12 .$ This model can be treated rigorously. Triangular noise was simulated as minimum of two uniform random numbers $u_1,u_2,$ and exponential noise as $1+\log u.$

\subsection{Pythagoras' theorem} 
As above for the covariance matrices, we now assign the six patterns 123, 132, 213, 231, 312 and 321 to the unit vectors in $\RR^6, $ that is $(1,0,0,0,0,0)$ up to $(0,0,0,0,0,1).$ Then the four pattern contrasts $\beta, \tau, \gamma$ and $\delta$ correspond to the vectors
\begin{eqnarray}\textstyle
\overline{\beta}=(1,0,0,0,0,-1)\ , \quad
\overline{\tau}=(\frac23,-\frac13,-\frac13,-\frac13,-\frac13,\frac23 ),  \notag
\\ 
\overline{\gamma}=(0,-1,1,1,-1,0),\quad \mbox{ and } \quad
\overline{\delta}=(0,1,1,-1,-1,0)\ .\label{vect}
\end{eqnarray}
Together with the vectors $c_1=(1,1,1,1,1,1)$ and $c_2=(0,-1,1,-1,1,0)$  of the two constraints  \eqref{cons1} and \eqref{cons2} these vectors form an orthogonal basis of $\RR^6.$  Their norms are
\begin{equation}\label{norms}
 |\overline{\beta}|^2=2,\quad |\overline{\tau}|^2=\frac43,\quad  |\overline{\gamma}|^2=|\overline{\delta}|^2=4\ . \end{equation}

Moreover, the four functions were defined so that they are all zero for white noise, and their value for pattern probabilities $p=(p_{123},...,p_{321})$ is $\tau(p)=\langle p,\overline{\tau} \rangle$ and similar for $\beta, \gamma, \delta .$ With $\langle x,y\rangle$ and $|x|$ we denote scalar product and Euclidean norm.

\begin{Proposition}[Pythagoras' theorem for pattern contrasts]\hfill
The distance to white noise defined in \eqref{d2} can be partitioned as follows:
\[ 4\Delta^2 = 3\tau^2 + 2\beta^2 +\gamma^2 + \delta^2 \ . \]
The quadratic distance of pattern frequencies $p$ to the pattern frequencies $p^{BM}$ of  symmetric random walk (see \eqref{brown3}) fulfils
\[ 4\cdot |p-p^{BM}|^2 = 3(\tau-\frac16)^2 + 2\beta^2 +\gamma^2 + \delta^2 \ . \]
For an arbitrary fixed pattern frequency distribution $p^o$ we have
\[ 4\cdot|p-p^o|^2 = 3(\tau(p)-\tau(p^o))^2 + 2(\beta(p)-\beta(p^o))^2 +(\gamma(p) -\gamma(p^o))^2 + (\delta(p)-\delta(p^o))^2 \ . \]
\end{Proposition}

\emph{Proof. } For an i.i.d.\  process all four pattern contrasts are zero, and for symmetric random walk only $\tau =1/6$ is non-zero. So the first two equations are special cases of the last one which we are going to prove now. The probability vectors $p$ and $p^o$ fulfil the constraints \eqref{cons1} and \eqref{cons2} so they are points in a four-dimensional affine space. Considering $p^o$ as origin of this space, we have a four-dimensional vector space with orthogonal basis $\overline{\tau},\overline{\beta},\overline{\gamma},\overline{\delta}.$ Now whenever $x$ is a vector in a four-dimensional space with an orthogonal basis $b_1,...,b_4,$ then 
\[  |x|^2 = \sum_{k=1}^4 \langle x, \frac{b_k}{|b_k|}\rangle^2 = \sum_{k=1}^4 \frac{1}{|b_k|^2} \langle x, b_k\rangle^2\ .  \]
Let $x=p-p^o$ and take $\overline{\tau},\overline{\beta},\overline{\gamma},\overline{\delta}$ as the $b_k.$ With \eqref{norms} we obtain the desired equation. \hfill$\Box$ \smallskip 

This statement is simple algebra, no statistics involved. It says that we can decompose the variance of the pattern distribution with respect to some origin or mean $p^o$ into components. We can say how much variance is due to $\tau , \beta$ etc. Thus we have new parameters which add to 1. In the case of an i.i.d.\  process  they are \cite{ba14}
\begin{equation}
 \tilde{\beta}=\frac{2\beta^2}{4\Delta^2},\quad \tilde{\tau}=\frac{3\tau^2}{4\Delta^2},\quad \tilde{\gamma}=\frac{\gamma^2}{4\Delta^2},\quad \tilde{\delta}=\frac{\delta^2}{4\Delta^2}\ . \label{relativ}
\end{equation}
This very much resembles the Anova methodology in statistics.  We now have to discuss the statistical independence of the components.

\subsection{Independence of estimators} \label{indi} 
We said that frequencies $p_\pi$ of different patterns $\pi$ are correlated, and the (asymptotic) covariance matrices near Brownian motion and white noise were calculated in \cite{els,sou,wei} and shown in \eqref{sigB} and \eqref{sigW}. 
Any two row vectors $x,y$ in $\RR^6$ can be interpreted as linear combinations of patterns, and $x\Sigma y'$ is the covariance of $x$ and $y.$ So matrix multiplication determines variances and covariances of our pattern contrasts $\beta,\tau, ...$ considered now as estimators over the space of all sample series of length $T.$ For \emph{symmetric random walk,} it turns out that all our four pattern contrasts are uncorrelated! Their variances are
\begin{equation}
T{\rm Var\,}\beta =1 ,\ T{\rm Var\,}\tau =\frac{1}{4} ,\ T{\rm Var\,}\gamma =\frac{1}{3},  \ T{\rm Var\,}\delta =\frac{2}{3}.
\label{variB}\end{equation} 
Principal component analysis of the covariance matrix $\Sigma_B$ is the standard method to find uncorrelated vectors. All eigenvalues of any covariance matrix are real and nonnegative. We determine the eigenvectors $x$ for each eigenvalue $\lambda .$ The variance of $x$ is $x\Sigma x'= \lambda xx'= \lambda |x|^2.$ And if $x,y$ belong to different eigenvalues then $x\Sigma y'=0, $ so they are uncorrelated.
In our case, the vectors $c_1,c_2$ of the two constraints \eqref{cons1} and \eqref{cons2} are eigenvectors corresponding to $\lambda=0,$ since they represent constant functions. 

\begin{Proposition}[Independence of pattern contrasts for symmetric random walk]\hfill
For the covariance matrix $T\cdot\Sigma_B$ of pattern frequencies of symmsetric random walk given in \eqref{sigB}, the vectors of the four pattern contrasts $\beta, \tau, \gamma, \delta$ are exactly the eigenvectors of the nonzero eigenvalues $\frac12, \frac{3}{16}, \frac16$ and $\frac{1}{12}.$ Their variances are given in \eqref{variB}.
\end{Proposition}
The proof is simple calculation, for instance $\overline{\beta}\Sigma_B=\frac12 \overline{\beta}.$ Note that eigenvectors are determined only up to a real constant, and eigenvector programs yield unit vectors. Our pattern contrasts are not unit vectors. Their variance is determined by $x\Sigma_Bx',$ or $\lambda |x|^2$ when we use  \eqref{norms}.

The case of an \emph{i.i.d.\  process} is a bit more complicated. By matrix multiplication with $\Sigma_W$ we get the (asymptotic) variances
\begin{equation}
T{\rm Var\,}\beta =\frac13 ,\ T{\rm Var\,}\tau =\frac{8}{45} ,\ T{\rm Var\,}\gamma =\frac{2}{5},  \ T{\rm Var\,}\delta =\frac{2}{3}.
\label{variW}\end{equation} 
The covariances are all zero, except for $T\cdot{\rm Cov\,} (\beta, \delta)=\overline{\beta}T\cdot\Sigma_W\overline{\delta}'=-1/3 .$ Thus 
\begin{equation}\label{corrbd}
{\rm corr\, } (\beta,\delta)=-\sqrt{2}/2=-0.707..
\end{equation}
Let us compare with eigenvectors.
In \cite{wei}, the nonzero eigenvalues of $T\cdot\Sigma_W$ are determined as  $\lambda_1=\frac{1}{12}(2+\sqrt{2}), \lambda_2=\frac{2}{15}, \ T\cdot\lambda_3=\frac{1}{10},$ and $\lambda_4=\frac{1}{12}(2-\sqrt{2}).$ Matlab yields the numeric eigenvalues $0, 0, 0.049, 0.1, 0.133, 0.285$ and the corresponding eigenvector columns

\begin{equation}\label{EV}
\left(\begin{array}{cccccc}
 0.12&0.39&-0.5&0&0.58&-0.5\\
-0.35&0.54&-0.35&-0.5&-0.29&0.35\\
0.60&0.24&-0.35&0.5&-0.29&0.35\\
-0.35&0.54&0.35&0.5&-0.29&-0.35\\
0.60&0.24&0.35&-0.5&-0.29&-0.35\\
0.12&0.39&0.5&0&0.58&0.5\\ 
\end{array}\right) \notag \end{equation}
 
The first two columns are linear combinations of the constraint vectors $c_1$ and $c_2.$ The columns 4 and 5 are the unit vectors of $\overline{\gamma}$ and $\overline{\tau},$ respectively. Thus $\gamma$ has variance $\frac{1}{10}\cdot 4=\frac25,$ and $\tau$ has variance $\frac{2}{15}\cdot \frac43=\frac{8}{45},$ as Bienaym\`e knew. The columns 3 and 6 are connected with $\beta$ and $\delta .$ Their sum is $-\overline{\beta},$ and their difference is $\overline{\delta}/\sqrt{2}.$ So the eigenvectors are $-\frac12\overline{\beta}\pm \frac12\overline{\delta}/\sqrt{2}.$ These are the orthogonal and uncorrelated unit vectors of $T\Sigma_W,$ and they are unique since the nonzero eigenvalues are all different.

It would be difficult to interpret these linear combinations, however.  We stick to $\beta$ and $\delta,$ which are represented by orthogonal vectors but have negative correlation  \eqref{corrbd}. (Another option would be to change $\delta$ for $\beta_2=\delta+\beta$ which is uncorrelated with $\beta$ and easy to interpret. But the vectors for $\beta$ and $\beta_2$ are not orthogonal, and Pythagoras' theorem then includes a mixed term.) 

\begin{Proposition}[Independence of pattern contrasts of an i.i.d.\  process]\hfill
For the covariance matrix $T\cdot\Sigma_W$ of pattern frequencies of white noise given in \eqref{sigW}, the vectors of pattern contrasts $\tau$ and $\gamma$ are eigenvectors of the eigenvalues $\frac{2}{15}$ and $\frac{1}{10}.$ The other two eigenvectors of nonzero eigenvalues generate the same plane as $\overline{\beta}$ and $\overline{\delta}.$ We choose $\beta$ and $\delta$ instead of eigenvectors, for better interpretation, and must take into account the strong negative correlation \eqref{corrbd}. All other pairs of the four estimators are uncorrelated.
\end{Proposition}

{\bf Open problem. } We have shown that for patterns of length 3, there is a canonical choice of four pattern contrasts. What about patterns of length 4?

\section{The descriptive power of turning rate}\label{V5}
\subsection{Order patterns and EEG} 
EEG brain signals are probably the widest field of application of order patterns. There is a well-developed methodology of EEG signals based on Fourier spectra and graphic elements, which is now combined with machine learning because of the size of datasets, see the comprehensive survey \cite{mach}. On the other hand, the data contain a lot of different artefacts which make spectral methods difficult. 
Order patterns apply most successfully to such really big and `dirty' datasets. 

Karsten Keller applied symbolic analysis as early as 2003 to EEG of epileptic children \cite{KL}, to assess the results of therapy performed by H. Lauffer. It turned out that larger permutation entropy indicates higher vigilance of the brain and thus a therapeutic success. 
He came back several times to such applications \cite{KUU,AKU}. Epilepsy is the most striking application of entropy methods in EEG \cite{RK,SL,sou} - there are too many papers to quote. Other studies concern the effect of anaestethic \cite{OSD} and psychotropic \cite{DMB} drugs, or Alzheimer's disease \cite{mor}. Permutation entropy and its variations are recommended methods because of their simplicity.  

Here we consider data from sleep research \cite{Te} which we had studied before \cite{ba17} with permutation entropy and $\Delta^2.$ It turns out that sleep analysis is still simpler with turning rate where we count only maxima and minima. This is a new application of pattern contrasts as function of time $t,$ not of the delay $d,$ in the setting of long data series. Turning rate will be used to compress the data, to classify sleep stages, and to find slow oscillations which are hardly accessible by conventional methods. 

\subsection{The dataset of Terzano et al.} 
We take the EEG sleep data of Terzano et al. \cite{Te} available on Physionet \cite{physio}. Although 20 years old, they are of excellent quality. There were several healthy volunteers which we shall study.  EEG were sampled with 500 Hz and a decent hardware low-pass filter at 30 Hz.  We shall take delay $d=4,$ that is 8 ms. In this range most EEG are totally disturbed by electrical noise from the power grid which is unavoidable in a clinical environment. It is common practice to filter away the `mains hum', which works well for conventional methods but can destroy the fine order structure. The Terzano data did not require any preprocessing.  

If one person is measured for 8 hours with 500 Hz, we have 15 million values. In each of 15 channels. For 15 persons. I have dealt with these data for several months, but cannot say I know them.  Such data must be evaluated automatically. Order patterns are not generally superior to machine learning but they have one advantage: we know what we do.

Nowadays high-resolutiion data are measured at every corner. Music is sampled with 40 kHz, laser experiments with Gigahertz frequency \cite{ba14}. Order patterns can screen the data, to see structure in high resolution. This is a chance to detect something new. 

\begin{figure}[h!] 
\begin{center}
\includegraphics[width=0.8\textwidth]{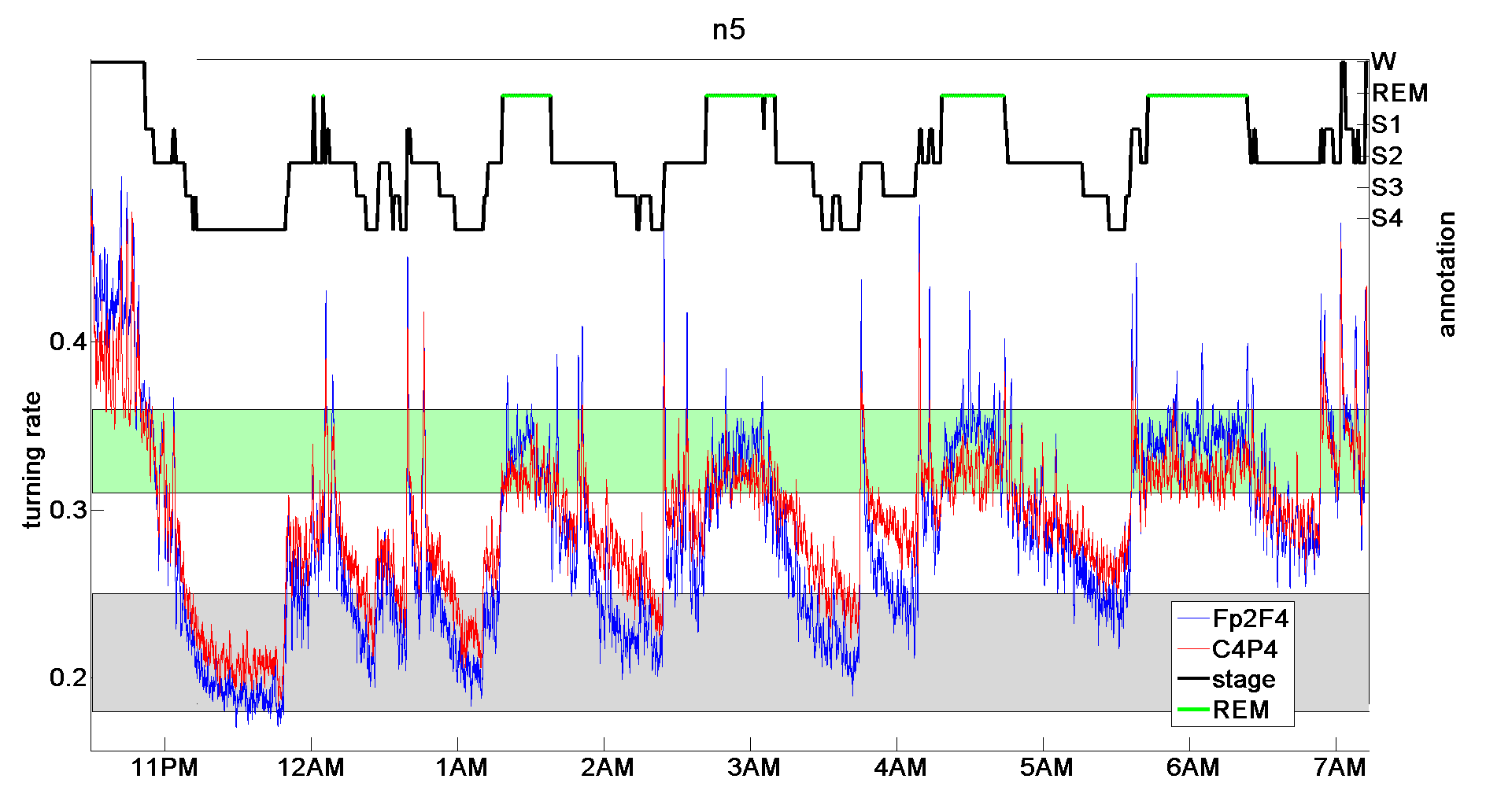}
\end{center}
\caption{The turning rate of an EEG directly classifies sleep stages}\label{TR1}
\end{figure}  

\subsection{Turning rate recognizes sleep stage}\label{sleepi}
Figure \ref{TR1} shows a sleep EEG of a healthy person. In the upper part you see the annotation of sleep stages by a medical expert, with one value for every epoch of 30 seconds.  Experts used various channels, among others eye trackers for REM phases. Their main goal was to find and analyze `cyclic alternating patterns'.  Annotation was an arts, with a comprehensive catalogue of rules. Now it is done by machine \cite{mach}, cheaper than experts.

\begin{figure}[h!] 
\begin{center}
\includegraphics[width=0.48\textwidth]{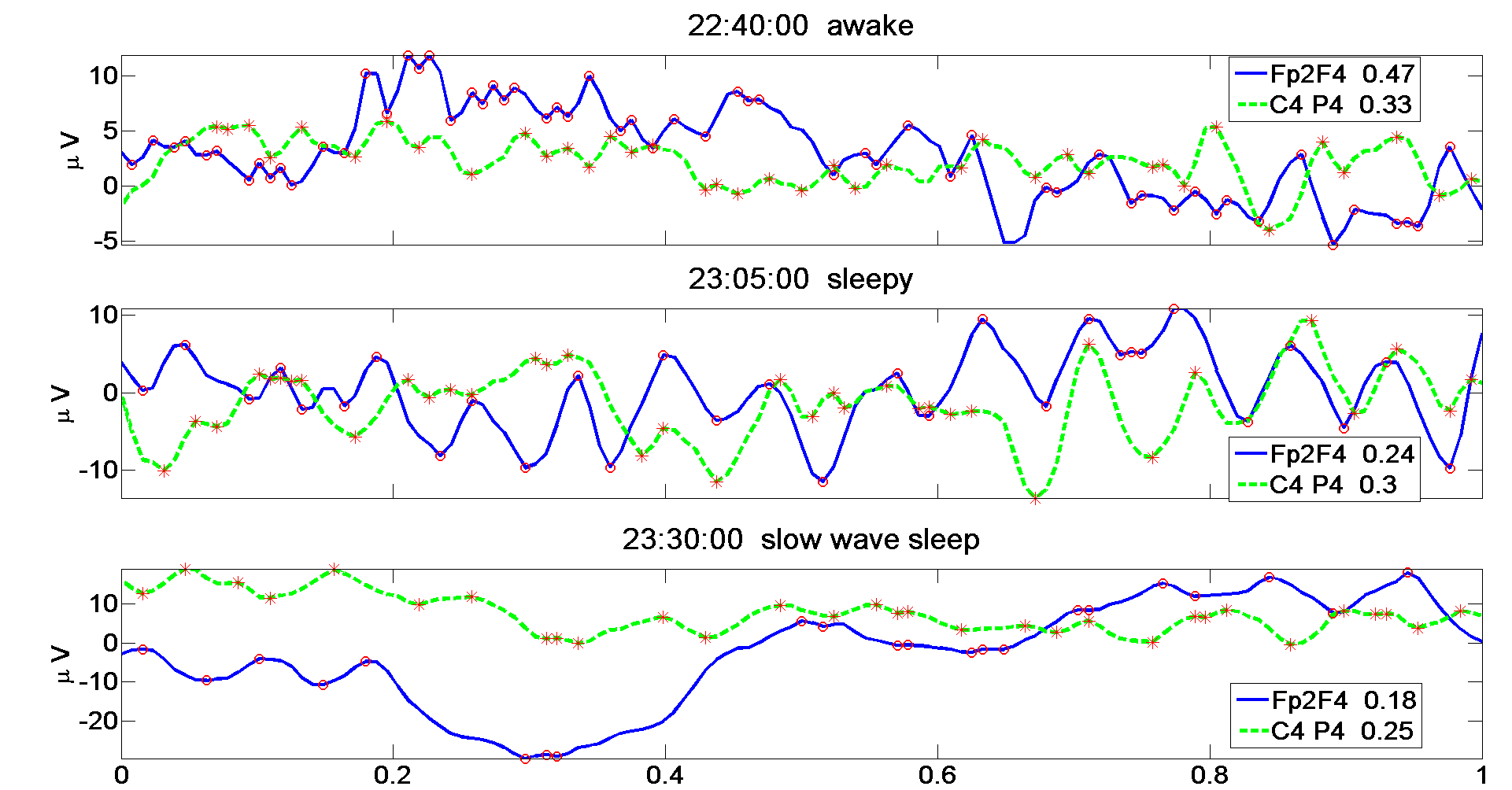}\quad
\includegraphics[width=0.48\textwidth]{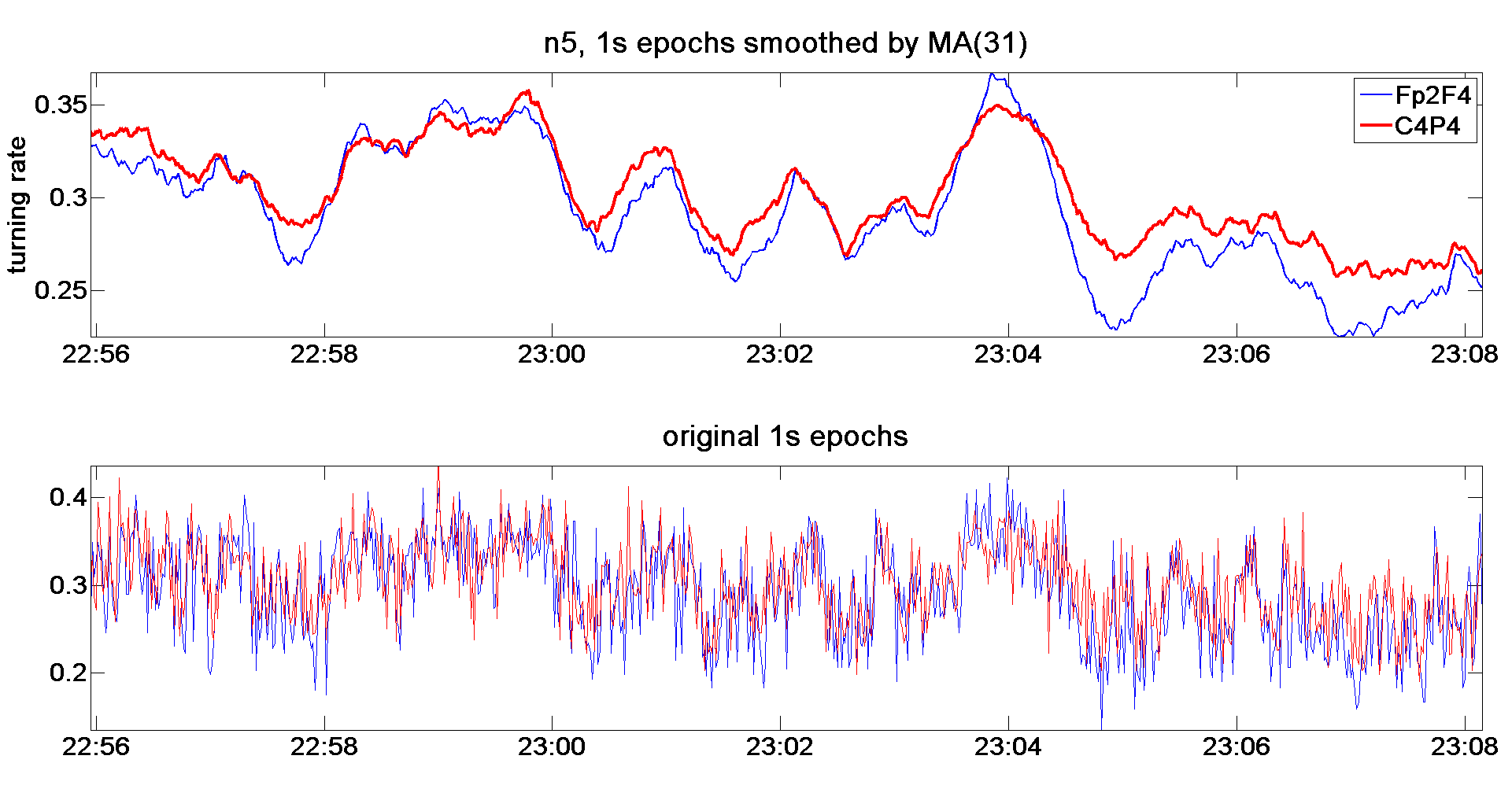}
\end{center}
\caption{Left: three clean data segments of one second demonstrate differences in turning rate. Right: turning rates of 1 s segments (bottom) are averaged over 30 s to find slow oscillations (top).}\label{TR2}
\end{figure}

In the lower part, two special EEG channel differences can be seen. Both functions are almost parallel to the annotation.  However, the two functions do not show the EEG values. They show the turning rate for $d=4.$ Thus by stupidly counting maxima and minima, we classify sleep almost as well as qualified experts!  In Figure \ref{TR2}, three very short parts of the time series demonstrate how the number of maxima and minima decreases in sleep.

The 15 million data of the very irregular original time series were compressed to a series of 1000 turning rates by analysing epochs of 30 s. With $T=15000$ values per epoch, the estimates of turning rate are fairly reliable even in the presence of small artefacts.  As a result, we obtain a new time series $\alpha_1,...,\alpha_{1000}$ of turning rates drawn in figure \ref{TR1}. No preprocessing, no filtering, no artefact removal.

\begin{figure}[h!] 
\begin{center}
\includegraphics[width=0.7\textwidth]{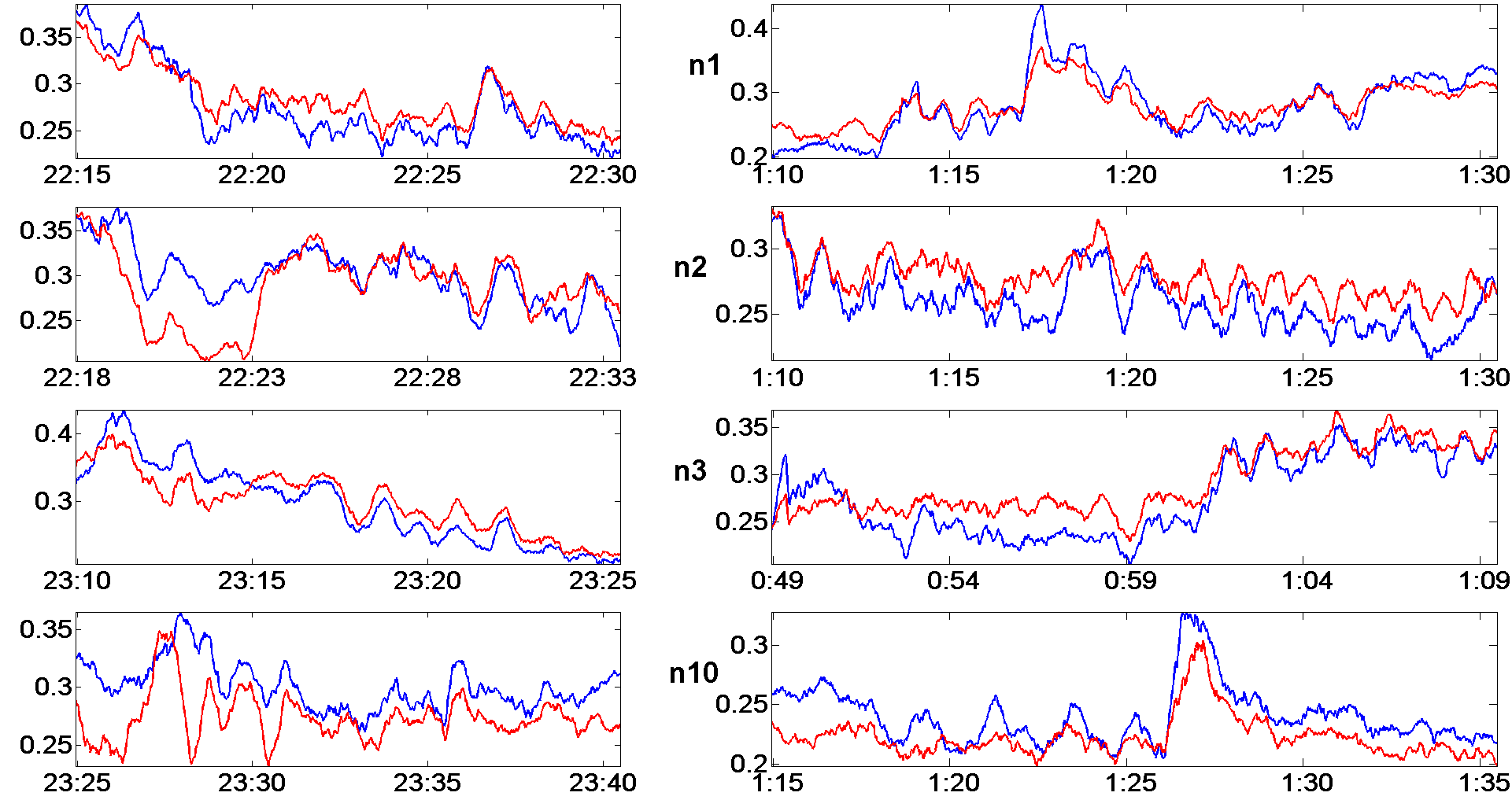}
\end{center}
\caption{Slow oscillations of turning rate for four other healthy subjects from the Terzano database. Note the synchrony between channels.}\label{ff15}
\end{figure}

Figures \ref{TR1}-\ref{ff15} are taken from an unpublished preprint \cite{ba18} where more than 10 persons from Terzano et al.\cite{Te} were studied, some with sleep disorders.  In spite of big differences between subjects, turning rate was in all cases tightly connected with sleep stage.  It seem to us that turning rate could even be used to {\em define sleep depth.} That would be better than a catalogue of diverse rules. Of course, this requires further study since $\alpha$ does not only depend on sleep stage, on $d$ and on the person, but also on measuring equipment, filtering options of hardware etc. \ REM phases of dream are more difficult to classify than sleep depth. Figure \ref{TR1} shows that the frontal channel Fp2F4 has higher turning rates in REM than the central channel C4P4 while outside REM it is the other way round. This effect was not seen for all persons.  
                      
\subsection{Infra-slow oscillations detected by turning rate}
So far, we took disjoint 30 s epochs to determine turning rates. For more detail, we shift the 30 s windows only by 1 s. Equivalently, we can calculate turning rates for 1 s segments and then apply a moving average with length 30, as shown in figure \ref{TR2}. The new time series has 30000 values $\alpha_t$ instead of 1000 and is fairly smooth.
 
We found clear oscillations of the turning rate with a period of about 1-2 minutes, at the onset of sleep, in light sleep (stage S2), and in the morning before getting awake.  Such infra-slow frequencies are totally beyond the scope of conventional Fourier analysis. See figure \ref{TR2} for our person n5, and figure \ref{ff15} for four other healthy persons. Both channels (and further channels not presented) show almost synchronous oscillations.
This may indicate a slow biorhythm from a central source and deserves further study. It is curious whether such rhythm is latently present in lecture rooms so that after one minute of concentration, students need a minute of rest.

\section{Some remarks on EEG data}

\subsection{Dependence of turning rate on the delay}
\begin{figure}[h!] 
\begin{center}
\includegraphics[width=0.48\textwidth]{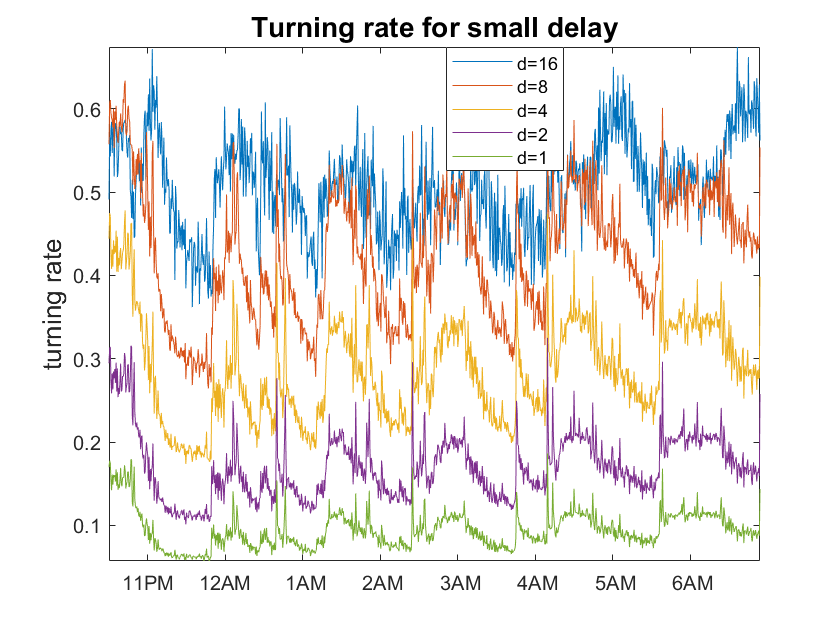}\quad
\includegraphics[width=0.48\textwidth]{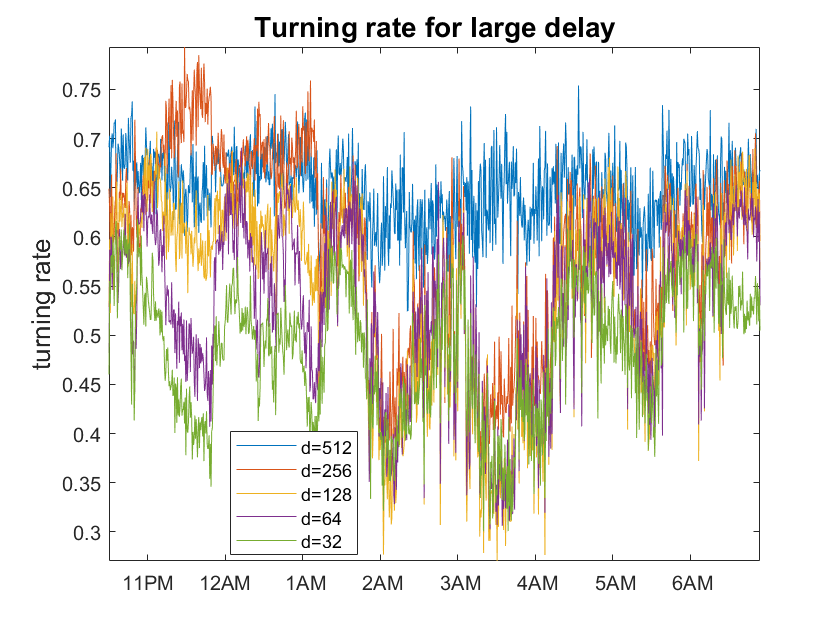}
\end{center}
\caption{Dependence of turning rates on the delay $d,$ ranging from 2 ms to 1 s.}\label{TR3}
\end{figure}

In the previous section we had fixed $d=4$ and considered $\alpha$ as a function of time. Now we ask how the $\alpha$-functions vary when we change $d.$
Figure \ref{TR3} illustrates the effects, for $d$ between 2 and 32 ms on the left, and for 64 ms up to 1 second on the right.  Actually, the $\alpha$ functions look very similar for most $d,$ and they would agree still better when we take $d=3,...,12,$ for instance. Thus we have a strong self-similarity in the order of EEG data.

For $d=1$ and $d=2$ the number of maxima and minima is very small. This comes from the hardware low-pass filter which smoothes the measured values, and will be similar for any physical measurement. The values of turning rate further increase with increasing $d,$ at least up to $d=32.$ This is quite natural. In figure \ref{TR2} we marked the maxima for $d=1.$ For $d=4$ there are much more maxima. 

With increasing $d$ the curves become more noisy. (This can be improved by taking windows according to `zero-crossing'.) The turning rate for $d=256$ and $d=512$ is near $2/3,$ the value for white noise. Here we have an effect of the high-pass filter which is also necessary for any physical measurement. It guarantees that the measured values do not move too far away from the zero line, as would be the case for random walk. The high-pass threshold was set to 0.3 Hz by Terzano et al..  When values move away for one second, the filter tends to get them back in the next second. In this way the self-similarity, present in models like Brownian motion, and apparently also in EEG, will be damped by filters on both sides in experimental data.

\subsection{Pattern contrasts in the EEG data}
After studying turning rate for fixed $d$ and sliding $t$-windows we now consider autocorrelation functions, with $d=1,...,500$ for fixed time epochs of 30 s. We now consider persistence $\tau(d)$ which resembles autocorrelation. Due to the low-pass filter, the time series is rather smooth on small scale and $\tau(1)$ is larger than 0.5, near to the maximal value $2/3.$ For white noise the value would be zero, for Brownian motion $\frac16 .$ Probably the first values of $\tau$ must be adjusted so that different measuring devices and hardware filters can be compared. There are no studies in this direction.

\begin{figure}[h!] 
\begin{center}
\includegraphics[width=0.9\textwidth]{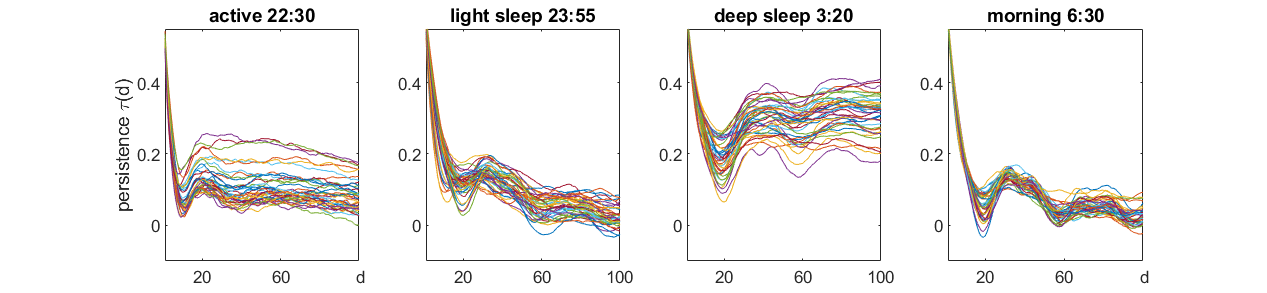}
\end{center}
\caption{Persistence in different stages. Each plot shows curves for 40 successive epochs of 30 s. The delay $d$ runs from 2 ms to 0.2 s.}\label{TR4}
\end{figure}

For small $d,$ persistence rapidly declines and reaches a minimum, in the first plot of Figure \ref{TR4} near $d_{\rm min}=10,$ in the other plots near $d_{\rm min}=20.$ This parameter is locked, the minimum was never at $d=15,$ say. The persistence at the minimum depends on the sleep stage and can be between 0 and 0.3. For $d\le 100,$ there were very few negative values of $\tau .$

Right of $d_{\rm min}$ the curve increases and reaches a maximum at $2d_{\rm min}.$ In many curves, a second minimum at $3d_{\rm min}\approx$ is visible. This indicates that we have here the so-called $\alpha$-rhythm, with frequency 12.5 Hz, corresponding to wave length 40 ms. This rhythm is always seen in occipital channels, and we had the frontal channel Fp2F4. Here it is apparent in the morning and did disappear most of the night, which is not shown in figure \ref{TR4}.  In deep sleep $\tau$ assumed large values up to $d=400.$ Thus we can classify this sleep stage also with large $d.$ We had chosen $d=4$ before the first minimum because the curves do not vary so much in this range.

Each of the plots in figure \ref{TR4} refers to 40 successive epochs, corresponding to 20 minutes of measurement. The $\tau$-curves seem to be rather consistent and interpretable. This was not true for the other contrasts $\beta, \gamma, \delta .$ Although they often show significant positive and negative values (see table \ref{testtab} below), the curves did not show consistent structure. They started near zero and looked like random walk. It seems that these contrasts express the numerous artefacts in EEG recordings. They may have some meaning, but not as average over windows of 15360 values. For instance, the cyclic alternating patterns which were the subject of the EEG study \cite{Te}, consist of an increasing oscillation and a decreasing oscillation. In our average of $\gamma$ over 30 s, the two will cancel out. 

We calculated the sum of squares of the contrasts over all epochs and $d=1,...,500,$ and determined their percentage according to \eqref{relativ}. We found that in the mean
\[ \tilde{\tau}=92\% \ ,\quad \tilde{\beta}=4.5\% \ ,\quad \tilde{\gamma}=1.3\% \ ,\ \mbox{ and } \tilde{\delta}=2.3\% \ .\]
Similar values were obtained when we restricted the domain of time and delays. Thus $\tau$ certainly plays the dominant rule.  To model EEG, one could start with assuming $\beta=\gamma=\delta=0.$

\subsection{Testing serial dependence of EEG}\label{serdep}
It is often necessary to test a time series against the null hypothesis that the values are from an i.i.d.\ process. In particular, any model of a process aims to explain the data in such a way that the residues, differences between observed and model values, cannot be distinguished from independent random numbers. As mentioned in section \ref{tailprob}, such tests on the base of order patterns were suggested by Matilla-Garcia and Ruiz Marin \cite{MM}, improved by Elsinger \cite{els}, and in this issue by Sousa and Hlinka \cite{sou} and Weiss \cite{wei}.  There is a great variety of other tests \cite{MM}. 

We shall not enter this matter, but discuss related questions for our EEG data: could they represent white noise, when they are sampled at some distance $d$? And will $H_3,H_4,$ and the pattern contrasts from section \ref{disti} distinguish them? We have fairly long time series with $T=15360$ for 1009 epochs of 30 s, and we can do the test for each $d=1,...,500$ between 2 ms and 1 s. Thus we can perform $500*1009=504500$ tests only for this dataset n5 from one night. Interdependence of tests and correction for multiple testing will be irrelevant for our questions. We are not looking for few tests which reject the i.i.d.\ null hypothesis. This should be the normal case with sample size $T=15360.$ We would like to have few tests which accept the null hypothesis. We use the 95\%\ confidence level to give all parameters a chance. Table \ref{testtab} shows that they all reject the i.i.d.\ process for the majority of cases, but there are big differences.

\begin{table}[h!]
\begin{tabular}{|c|c|c|c|c|c|c|c|c|c|}\hline
Test parameter &\multicolumn{2}{c|}{$\log{24}-H_4$}&$\log {6}-H_3$&$\tau$&$\beta$&$\gamma$&$\delta$&$\tau_4$&$\beta_4$\\ \hline                   
confidence level \%&95&99.9&95&95&95&95&95&95&95\\ \hline
$H_0$ accepted \%&0&0.002&0.5&6.7&24&36&32&9.2&22 \\  \hline
significantly larger&100&99.998&99.5&80&30&40&37&73&30 \\ \hline
significantly smaller&-&-&-&14&46&24&31&17&49\\   \hline
\end{tabular} 
\caption{Testing EEG data with an i.i.d.\ null hypothesis} \label{testtab}
\end{table}

For permutation entropy we took the critical quantiles for 95\%  and 99.9\%  from table \ref{tailq}. For pattern contrasts we took standard normal quantiles, multiplied by $\sigma/T$  where $\sigma^2$ is taken from \eqref{variW}. Compare Weiss \cite{wei} where simulations justified the normal approximation already for much smaller $T.$ We also included the new tests of Weiss, Ruiz Marin, Keller and Matilla-Garcia \cite{MW} with two pattern contrasts of length 4 as follows. Let $\beta_4=p_{1234}-p_{4321}$ and $\tau_4=p_{1234}+p_{4321}-\frac{1}{12}.$ In \cite{MW} it is shown that for an i.i.d.\ process $(p_{1234},p_{4321})$ has an asymptotic normal distribution with covariance matrix $\frac{1}{4032}{199 \ -17\choose -17 \ 199}.$ 
Thus under the null hypothesis $\beta_4$ and $\tau_4$ have mean zero and variances $\frac{3}{28}$ and $\frac{181}{2016},$ respectively.

Table \ref{testtab} shows that permutation entropy has much better power than the contrasts in these serial dependence tests. Even $H_3$ is much better than $\tau$ although simulations of Weiss \cite{wei} indicate that $\tau$ can perform better for small $T.$ It seems that for large $T,$ the $1/T$-scaling of $H_m$ becomes superior to the $1/\sqrt{T}$-scaling in the normal distribution. For our EEG data, all pattern contrasts show significant deviations to both sides from an i.i.d.\ process. Of course this can change when we consider special conditions, for instance deep sleep, where $\tau$ and $\tau_4$ are large (small turning rate, see section \ref{sleepi}).  In this connection, we note that for $d=1,...,50,$ all 50450 tests with $\tau_4$ gave significant positive deviation from the null hypothesis.  This may be due to the influence of the hardware low-pass filter, however. Also $\tau$ accepted $H_0$ in less than 0.1\%  of tests with $d\le 50.$

These results support the original suggestion of Matilla-Garcia and Ruiz Marin \cite{MM} to use $H_4$ as test quantity. By definition, it can detect all differences in pattern frequencies for $m\le 4.$ On the confidence level of 99.9\% , there were only 9 from 504500 tests where $H_0$ was accepted. This means: 
there is not the slightest chance to find an i.i.d.\ process in the EEG. These data represent another type of noise, and new parametric models are needed to understand them.

Of course, when monthly economic data are given, we have short time series with $T\le 300$ and can just estimate patterns of length 3. Use of an i.i.d.\  process is justified in such a context. But for high-resolution measurements, we should work with patterns of length 4 and should develop better models for noise so that we can dismiss the i.i.d.\  process.  However, fractional Brownian noise is too special to provide a model for EEG. And $1/f$-noise \cite{Man} is a physical phenomenon which actually includes our EEG, but not a model. So far, diverse mathematical explanations of $f^{-\beta}$-noise have not found acceptance \cite{scholar}. In the next sections, we try another way.

\section{A new kind of ordinal process}\label{V6}
\subsection{Ranking fighters by their strength.}  
So far, order patterns have been derived from numerical values of a time series. Now we go the opposite way. This occurs in the world of tennis or chess were players first pairwise compare each other and afterwards are assigned a rank or ELO score.

Suppose there are four players $A, B, C, D$ in a tennis club. We design an algorithm to order them by strength. First $B$ plays against $A.$ Let us assume $B$ wins. Then we write $A<B.$  Now the next player $C$ plays against the previous player $B.$  If $C$ wins, then $A<B<C$ by transitivity of order. If $C$ loses, however, $C$ must still play against $A.$ Depending on the result, either $A<C<B$ or $C<A<B.$ Now the last player $D$ comes and plays first against the previous player $C,$ and then against $B$ and/or $A,$ if the comparison is not yet determined by previous games and transitivity.

This procedure is now turned into a mathematical model. We want to construct a random order between objects $X_1, X_2,...$ which are not numbers.  We follow the algorithm for tennis players, and replace each match by the toss of a fair coin. As a result we obtain a stochastic process on ordinal level, like white noise or Brownian motion. It will be possible to determine pattern probabilities much better than for Brownian motion.

\subsection{The coin tossing order.} \label{coin} 
We define a random order between objects $X_1, X_2, X_3,...$ by tossing fair coins. We write 1 for `head' and 0 for `tail'.  Thus our basic probability space $\Omega$ is the space of all 0-1-sequences, where each coordinate is 0 or 1 with probability $\frac12 ,$ independently of all other coordinates.  We write
\[ \Omega = \{ (c_{21}, c_{32}, c_{31}, c_{43}, c_{42},...)\, |\, c_{ji}\in\{ 0,1\} \} \]
where $c_{ji}=0$ means $X_i<X_j,$ and $c_{ji}=1$ means $X_i>X_j,$ for $1\le i< j.$ The important point is that $c_{ji}$ will be disregarded when the order of $X_i$ and $X_j$ is already fixed by previous comparisons and transitivity.

The first coin $c_{21}$ decides the ordering of $X_1$ and $X_2.$ Now suppose $X_1,...,X_{j-1}$ are already ordered. Then $X_j$ is compared to $X_{j-1}, X_{j-2},..., X_1$
by considering the random numbers $c_{ji}.$ However, when the comparison is fixed by transitivity from the already defined ordering, then $c_{ji}$ is disregarded - that coin would not be thrown.  

\begin{figure}[h] 
\begin{center}
\includegraphics[width=0.7\textwidth]{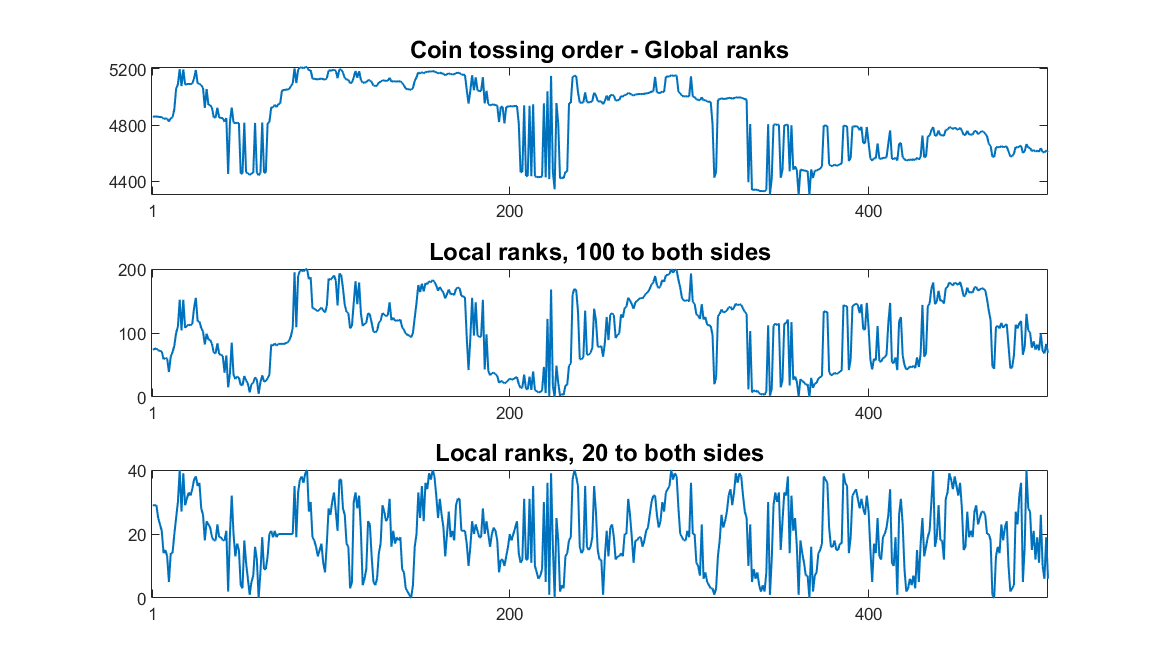}
\end{center}
\caption{Global and local rank numbers obtained from coin tossing.}\label{CTOfig}
\end{figure}  

The resulting random order will  be called \emph{coin tossing order.} It can be easily simulated. Figure \ref{CTOfig} shows rank numbers of 500 consecutive objects $X_j$ in th middle of a simulated series  of length $T=10000.$ Global rank numbers have strange discontinuities but local rank numbers obtained by comparing with the next 20 objects on the left and right show a familiar picture.

\subsection{Properties of the coin tossing order.} 
\begin{Proposition}[Basic properties of the coin tossing order]\hfill \label{coip}\begin{enumerate}
\item[ (i)] For any permutation $\pi\in S_m,$ the probability that $X_t,...,X_{t+m-1}$ shows pattern $\pi$ is of the form $2^{-u}$ for some integer $u=u(\pi).$
\item[ (ii)] $P\{ X_i<X_j\}=\frac12$ for any time points $i<j.$ 
\item[(iii)] The coin tossing order is order stationary and has the Markov property.
\item[(iv)] The pattern probabilties $P_\pi (d)$ for patterns of length 3 are the same for Brownian motion and coin tossing order. Thus $\beta(d)=\gamma(d)=\delta(d)=0$ and $\tau(d)=\frac16 ,$ for all $d\ge 1.$
\end{enumerate}
\end{Proposition}

\emph{Proof. } (i): $u(\pi)$ is the number of coin flips needed to determine $\pi .$ \quad 
(ii): The transformation $\varphi :\Omega\to\Omega, \varphi(...c_k...)=(...(1-c_k)...),$ which interchanges `head' with `tail' for all coins, preserves the probability. Moreover, $\varphi (\{ X_i<X_j\} )= \{ X_i>X_j\} ,$ so the two events have the same probability. \quad
(iii): The probability $P_\pi^t$ that pattern $\pi$ appears at $X_t,...,X_{t+m-1}$ is determined by the coins 
$c_{ji}$ with $t\le i<j\le t+m-1$ according to rules that do not depend on $t.$ So $P_\pi^t$ does not depend on $t.$ This holds for $d=1$ and all $m>1,$ which implies order stationarity also for $d>1.$ Moreover, when $t+m-1\le s$ then $P_\pi^t$ is determined from the $c_{ji}$ with $j\le s,$ and $P_\pi^s$ is determined by the  $c_{ji}$ with $j>s.$ So the occurence of patterns up to $s$ and after $s$ is independent. This is the Markov property.\quad
(iv): In virtue of the last two properties, the calculation \eqref{brown3} holds true for coin tossing order, for any $d.$
\hfill$\Box$ \smallskip 

\begin{figure}[h!] 
\begin{center}
\includegraphics[width=0.48\textwidth]{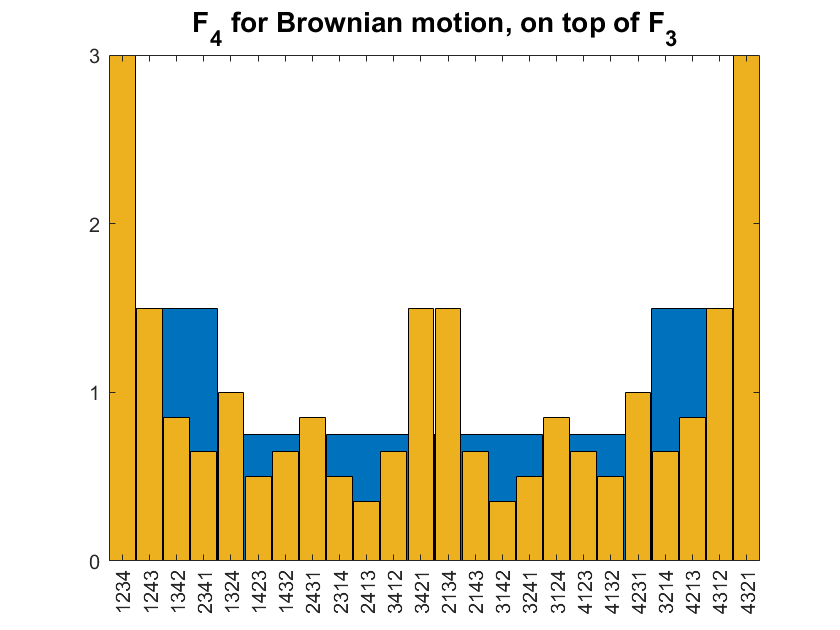}\quad\includegraphics[width=0.48\textwidth]{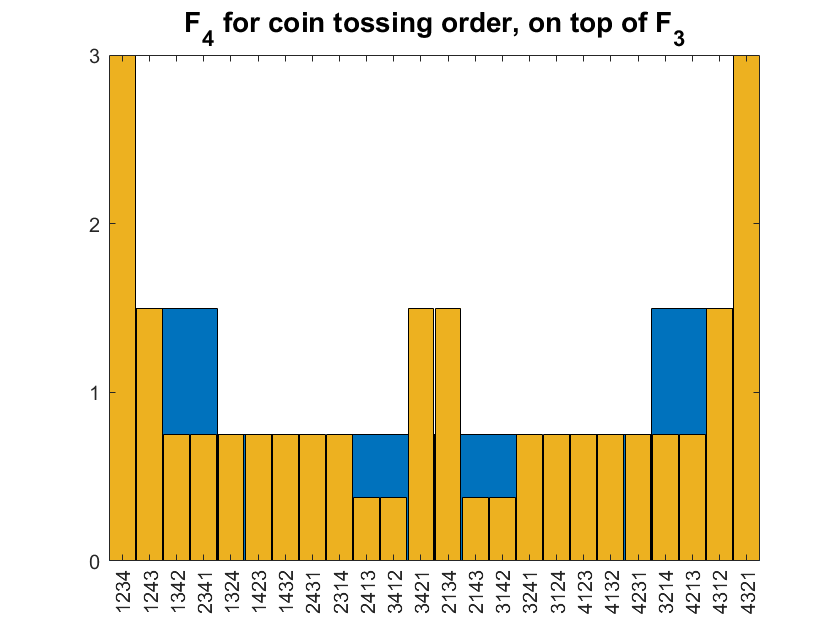}
\end{center}
\caption{Probabilities of patterns of length 4 (brown) on top of those of length 3 (blue) for Brownian motion (left) and coin tossing order (right). For length 3, the probabilties of the two processes coincide.}\label{bc4}
\end{figure}  

As a consequence of (ii) and (iii), the probability of the increasing pattern $12...m$ is $2^{1-m}$ for any $m$ and both Brownian motion and coin tossing order.  However, the distributions of patterns of length 4  are quite different, as can be seen in figure \ref{bc4}. Coin tossing looks a bit artificial and is indeed not a good model for financial time series, cf. \cite{ba19}. This figure shows that we must investigate patterns of length $\ge 4$ if we want to distinguish processes like these two.  Moreover, Brownian motion is self-similar so that the probabilities $P_\pi (d)$ do not depend on $d.$  Computer calculation for coin tossing order showed that there is no self-similarity. For $\pi =1324$ we have $P_\pi (1)=\frac{1}{32}\approx 3.12\% $ while $P_\pi (2)\approx 3.42\% $ and $P_\pi (3)\approx 3.49\% .$ 

{\bf Open problem. } Let us define new pattern probabilities $Q_\pi (1)=P_\pi (d)$ for $d=2,3,...$ What are the properties of $Q$ ?  Will there be a limit for $d\to\infty$ ?\smallskip

The function $u(\pi)$ in (i) can be expressed in closed form: it is the number of all $(i,j)$ with $1\le i<j\le m$ such that there exists no $k$ with 
\[ i<k<j \quad\mbox{ and }  \quad \pi_i<\pi_k<\pi_j \ \mbox{ or } \ \pi_j<\pi_k<\pi_i\ .\]
This allows to determine easily all pattern probabilities for $d=1$ and $m\le 10$ by computer.  It also shows that coin tossing order is reversible, that is $\pi_1\pi_2...\pi_m$ and $\pi_m...\pi_2\pi_1$ have the same probability.
Furthermore,  $u(\pi)$ can be interpreted as an energy function. Let $Q$ denote the probability measure on $S_m$ given by the probabilities $p_\pi$ of coin tossing order, and let $M_Q(u)$ denote the mean energy with respect to $Q.$ Then 
\[  p_\pi = e^{-u(\pi)\cdot\log 2} \quad \mbox{ and }\quad  H(Q)=M_Q(u)\cdot\log 2 \ .\]
For all other probability measures $Q'$ on $S_m,$ the permutation entropy is smaller than the mean energy. Thus $Q$ is a so-called Gibbs measure on $S_m,$ see chapter 1 in  \cite{gk}. 

{\bf Open problem. } Are there other meaningful energy functions for permutations, perhaps even parametric families of Gibbs measures on $S_m$ ?

\section{Rudiments of a theory of ordinal processes}\label{V7}
\subsection{Random stationary order}\label{last} 
Having studied one example, we turn to discuss a possible theory of ordinal stochastic processes. Now $X_1,X_2,...$ will not be numbers, only objects which are ordered. We still have the discrete time domain $\{ 1,...,T\}$ for time series and $\NN =\{1,2,...\}$ for models. An order will be a relation $<$ on the time domain such that $x\not< x,$ and $x<y$ implies $y\not<x,$ and $x<y$ plus $y<z$ implies $x<z$. Of course, the order does not apply to the time points $t,$ but to the corresponding objects $x_t.$ Then the property that $x_t,...,x_{t+m-1}$ shows pattern $\pi$ has a meaning for every $\pi\in S_m$ -- it is either true or false.  When an order on $\NN$ or $\{ 1,...,T\}$  is given, we can determine pattern frequencies, permutation entropy and so on. In the following, we always have $d=1.$ 

We want to construct models of ordinal processes, like the coin tossing algorithm. For this purpose we need the concept of random order. To keep things simple, \emph{a random order is defined as a probability measure on the set of all orders on the time domain.} For the finite time domain $\{ 1,...,T\} ,$ a random order is just a probability measure $P_T$ on the set $S_T$ of permutations of length $T.$ For $m\le T$ and $\pi\in S_m$ and $1\le t\le T+1-m$ the random order allows to determine the probability
\begin{equation}\label{ptpi} 
P_\pi^t = P\{ x_t,...,x_{t+m-1}\mbox{ shows pattern } \pi \} 
\end{equation}
The random order will be called \emph{stationary} if the  $P_\pi^t$ do not depend on the time point $t,$ for any pattern $\pi$ of any length $m<T.$ In other words, the numbers $P_\pi^t$ must be the same for all admissible $t.$ This is exactly the order stationarity which we defined for numerical processes in \eqref{ost}.

\subsection{The problem to find good models}
The infinite time domain $\NN$ is considered in the next section. Using classical measure theory, that will be easy. The real problem appears already for finite $T,$ even for $T=20.$ We have an abundance of probability measures on $S_T$ since we can define $P_\pi$ for every $\pi .$ When we require stationarity, we shall have $(T-1)!$ equations for the $T!$ parameters, as shown below, which is still too much choice.  

The problem is to select realistic $P_\pi .$ Most of the patterns $\pi$ for $m=20$ will never appear in any real time series, and we could set $P_\pi=0.$  But we do not know for which $\pi .$ There are three types of properties which we can require for our model.

\begin{itemize}
\item Independence: Markov property or $k$-dependence of patterns (cf. \cite{sou}).
\item Self-similarity: $P_\pi(d)$ should not depend on $d.$
\item Smoothness: There should not be too much zigzag in the time series. Patterns like 3142  should be exceptions.
\end{itemize}

All three properties contribute to the main goal - a simple, comprehensible model definition. An algorithm like coin tossing is not necessary, but it may help. We only need the $P_\pi .$

{\bf Open problem. } Does there exist on $S_{100}$ a stationary random order which is Markov and self-similar (for admissible $d$ and $\pi$) and has parameters different from  white noise and Brownian motion?  For instance $\beta\not= 0,$ or $0<\tau<\frac16$ ?

\subsection{Random order on $\NN$}
The infinite time domain $\NN$ has its merits. Order stationarity, for instance, is very easy to define since every pattern can be shifted to the right as far as we want. It is enough to require
\begin{equation}\label{statinf}
P_\pi^t = P_\pi^{t+1} \quad\mbox{ for all finite patterns }\pi \mbox{ and } t=1,2,...
\end{equation}
It is even enough to require $P_\eta^1 = P_\eta^2$ for all patterns $\eta$ of any length. (To prove that this implies \eqref{statinf} for a fixed $t$ and a pattern $\pi$ of length $k,$  consider all patterns $\eta$ of length $m=t+k-1$ which show the pattern $\pi$ on their last $k$ positions. The sum of their probabilities $P_\eta^1$ equals $P_\pi^t$ since $P$ is a measure.  And shifting all $\eta$ from 1 to 2 means shifting $\pi$ from $t$ to $t+1.$)

On the other hand, infinite patterns require a limit $T\to\infty .$
Actually, there are lots of recent papers on infinite permutations, that is, one-to-one mappings of $\NN$ onto itself.  An overview is given in Pitman and Tang \cite{PT}. However, an order on $\NN$ is a much wider concept than a permutation on $\NN .$ An infinite permutation $\pi_1,\pi_2,...$ defines an order on $\NN ,$ with the special property that below a given value $\pi_k$ there are only finitely many other values, for any $k.$ 

For an order on $\NN ,$ however, there need not exist smallest objects - usually each object has infinitely many other objects below and above.
Nevertheless, an order on $\NN$ is uniquely defined by patterns $\pi^m\in S_m$ which represent the order of the first $m$ elements, for $m=2,3,...$  For example
$ 1\, 2 - 2\, 3\, 1 - 3\, 4\, 1\, 2 - 4\, 5\, 1\, 3\, 2 -...$

\begin{Proposition}[Approximating random order on $\NN$]\hfill \label{app}\begin{enumerate}
\item[ (i)] A sequence of permutations $\pi^m\in S_m$   defines an order on $\NN$ if for $m=2,3,...,$ the pattern $\pi^m$ is represented by the first $m$ values of the pattern $\pi^{m+1}.$ 
\item[ (ii)] A sequence of probability measures $P^m$ on $S_m$   defines a random order $P$ on $\NN$ if for $m=2,3,...,$ for $\pi\in S_m$ and $\pi^{m+1}\in S_{m+1}$ holds
\[ P_{m+1}\{ \pi^{m+1} \mbox{ shows pattern $\pi$ at its first $m$ positions }\} =P_m(\pi) \ .\]
\item[(iii)] The random order $P$ defined by the $P_m$ is stationary if and only if for $m=2,3,...,$ for $\pi\in S_m$ and $\pi^{m+1}\in S_{m+1}$ holds
\[ P_{m+1}\{ \pi^{m+1} \mbox{ shows pattern $\pi$ at its last $m$ positions }\} =P_m(\pi) \ .\]
\end{enumerate}
\end{Proposition}

\emph{Proof. } (i): The condition says that the pattern of the first $m$ objects, defined in step $m,$ will not change during successive steps. So the construction is straightforward. The rank numbers, however, may change in each step, and they may converge to $\infty$ for $m\to\infty .$ \quad 
(ii): The condition says that the probability for a pattern $\pi\in S_m$ to appear for the first $m$ objects, defined by $P_m,$ will remain the same for $P_{m+1}$ and successive probability measures. So we can define 
\[ P\{ \mbox{ the objects with numbers $1,...,m$ show pattern  $\pi$ }\} = P_m(\pi )\]   
in a consistent way, and the $P_m$ determine $P.$ Actually, $P$ is an inverse limit of the measures $P_m.$ Below, we provide a more elementary argument.
\quad             
(iii): Together with (ii), this condition says that $P_\pi^1 = P_\pi^2$ for all patterns $\pi$ of any length $m.$ As noted above, this is just the definition \eqref{statinf} of stationarity.
\hfill$\Box$ \smallskip 

\subsection{The space of random orders.}
We show that the set of all orders on $\NN$ can be represented by the unit interval, using a numeration system similar to our decimal numbers. We first assign subintervals of $I=[0,1]$ to the permutations of length $2,3,...$  The pattern 12 will correspond to $[0,\frac12],$ and 21 to $[\frac12,1].$ The permutations 123, 132, and 231, which show the pattern 12 at their first two places, will correspond to $[0,\frac16], [\frac16,\frac13],$ and $[\frac13,\frac12],$ respectively. For intervals of length 4, see figure \ref{bc4}.

Instead of lexicographic order, we define a hierarchical order of permutations, with respect to the patterns shown by the first 2,3,... elements. For any permutation $\pi=\pi_1...\pi_m$ of length $m$ and any $k$ with $2\le k\le m$, let $r_k(\pi)$ denote the number of $j<k$ with $\pi_j>\pi_k.$
This is a kind of rank number of $\pi_k$ with values between 0 and $k-1.$ For example $r_3(123)=0, r_3(132)=1,$ and $r_3(231)=2$ while $r_2=0$ for these three patterns. 
Now we assign to the permutation $\pi$ of length $m$ the following interval:
\begin{equation}\label{hier}
I(\pi)=[x(\pi), x(\pi)+\frac{1}{m!}] \quad \mbox{ with } \quad  x(\pi)=\sum_{k=2}^m \frac{r_k}{k!} \ . \end{equation}
It is easy to check that the patterns $\pi^{(k)}=\pi_1...\pi_k$ of the first $k<m$ items of $\pi$ are assigned to larger intervals:
\[  I\supset I(\pi^{(2)}) \supset I(\pi^{(3)}) \supset ...\supset I(\pi^{(m)}) \]
where $\pi^{(m)}=\pi .$ Actually, $\pi$ need not be a permutation, just a pattern - it could also be a numerical time series. Only the ordering of the $\pi_j$ is used for defining $r_k$   and $I(\pi).$ 

When we extend the pattern to the right, we get smaller nested subintervals, and for $m\to\infty$ a single point $x=\sum_{k=2}^\infty \frac{r_k}{k!}$ which characterizes the limiting order of infinitely many objects. Thus each order on $\NN$ corresponds to a unique point $x$ in $[0,1].$  This is very similar to decimal expansions where we subdivide an interval into 10 subintervals. In case of patterns, the $r_k$ are the digits, and we subdivide first into 2, then 3, 4, 5,... intervals. The endpoints of intervals represent two orders on $\NN ,$ but this is an exception, as 0.5 and $0.4999...$ for the decimals.

Once we have represented all orders on $\NN$ as points in $[0,1],$ we can better understand the probability measures $P_2, P_3, ...$ of Proposition \ref{app} and the limiting probability measure $P$ which is called random order on $\NN .$ We start with the function $F_1(x)=1$ which denotes uniform distribution on $[0,1].$  The function $F_m$ will represent the measure $P_m,$ for $m=2,3,...$  For patterns $\pi$ of length $m$ it is defined as histogram of $P_m:$
\begin{equation}\label{pm} 
F_m(x)=m!\cdot P_m(\pi) \quad\mbox{ for }    x\in I(\pi) \ . \end{equation}
See figure \ref{bc4}. The rectangle over $I(\pi)$ has area $P_m(\pi).$ In case of white noise, $F_m=1$ for all $m,$ and the $\lim F_m=F=1$ is the uniform distribution. We now show that such limit exists for all sequences $F_2,F_3,...$ for which the $P_2,P_3,...$ fulfil condition (ii) of Proposition \ref{app}. We reformulate (ii) as 
\begin{equation} \label{mar}
\int_{I(\pi)} F_{m+1}(x) dx = m!\cdot P_m(\pi)  =\int_{I(\pi)} F_m(x) dx \ . \end{equation}
The second equation is obvious, and the first is best shown by example, for $m=3$ and $\pi=312.$ We have $r_2=r_3=1,$ so $I(\pi)=[\frac46,\frac56 ].$ The possible extensions $\pi^{4}$ of $\pi$ are 3124, 4123, 4132, and 4231. Their intervals of length $\frac{1}{24}$ partition $I(\pi).$ Condition (ii) says that $P_4\{ 3124, 4123, 4132, 4231\} =P_3\{ 312\} .$
Thus the four rectangles of $F_4$ over $I(\pi)$ together have the same area as the one rectangle of $F_3.$ This is expressed in \eqref{mar}.

Equation \eqref{mar} says that the $F_m$ form a martingale. The martingale convergence theorem implies that there is a limit function $F,$ in the sense that $\int_0^1 |F_m-F| dx$ converges to zero for $m\to\infty .$ This limit function is integrable and $\int F=1.$ As a density function, it defines the probability measure $P$ on all orders on $\NN .$

Our argument indicates that random orders on $\NN$ belong to the realm of classical analysis and probability.  Of course, the density function $F$ will be terribly discontinuous and can hardly be used to discuss stationarity.  (Open problem for experts: find examples for which the $F_m$ converge in $L_2.$)

\subsection{Extension of pattern distributions} 
We conclude our paper with an optimistic outlook for practitioners. When you find a distribution of pattern probabilties of length 3 or 4, from data say, you need not care for an extension to longer finite and infinite patterns. Such an extension will always exist.
  
\begin{Proposition}[Markov extension of pattern probabilities]\hfill \label{ext}
Any stationary probability measure $P_m$ on $S_m$ can be extended to a stationary probability measure $P_{m+1}$ on $S_{m+1}$, and hence also to a stationary probability measure $P$ on the space of random orders.
\end{Proposition}

\emph{Proof. } 
To show that  $P_{m+1}$  is stationary, we need only   verify condition (iii) of Proposition \ref{app}: Any pattern within $\{ 1,...,m\}$ can be shifted to the right, without changing its probability, until its maximum reaches $m,$ by assumption of stationarity of $P_m.$ It remains to shift the maximum from $m$ to $m+1,$ and this is condition (iii).

In the following extension formula, we use the convention that probability measures on permutations also apply to patterns, by replacing a pattern with its representing permutation.
Let $\pi=\pi_1\pi_2...\pi_m\pi_{m+1}$ be a permutation in $S_{m+1}.$
\begin{equation}\label{ex}
P_{m+1}(\pi)=\frac{P_m(\pi_1...\pi_m)\cdot P_m(\pi_2...\pi_{m+1})}{P_m(\pi_2...\pi_m)}
\end{equation}
This formula will be used whenever there exists some $\pi_k$ with $2\le k\le m$ between $\pi_1$ and   $\pi_{m+1}.$ However, if $\pi_1$ and $\pi_{m+1}$ are neighboring numbers, then the right-hand side of \eqref{ex} cannot distinguish $\pi$ and $\pi '=\pi_{m+1}\pi_2...\pi_m\pi_1.$ In such cases, both $\pi$ and $\pi '$ are assigned half of the value of the right-hand side, in order to avoid double counting.

The denominator on the right refers to a pattern $\pi_2...\pi_m$ which has a representing permutation $\kappa=\kappa_1...\kappa_{m-1}$ in $S_{m-1}$ and can be extended in $m$ ways to a permutation $\eta=\eta_1...\eta_m \in S_m.$ Indeed, $\eta_m$ can be chosen from $\{ 1,...,m\} .$ For $j<m$ then either $\kappa_j<\eta_m$ and  $\eta_j=\kappa_j,$ or 
$\kappa_j\ge \eta_m$ and $\eta_j=\kappa_j +1.$ Let us write
\[ P_m(\kappa)=P_m\{ \eta\, |\, \eta_1...\eta_{m-1} \mbox{ shows pattern } \kappa \} 
\sum_\eta P_m(\eta) \ . \]
In the numerator, $\pi_1...\pi_m$ is also a pattern which has a representing permutation $\lambda=\lambda_1...\lambda_m$ in $S_m.$ This term is $P_m(\lambda).$

We now prove that the defined $P_{m+1}$ fulfils condition (ii) of Proposition \ref{app}. We calculate
\[ p=P_{m+1}\{ \pi_1...\pi_{m+1}\, |\,  \pi_1...\pi_m \mbox{ is represented by } \lambda \} \]
using the definition \eqref{ex} of $P_{m+1}.$ There are $m+1$ permutations $\pi\in S_{m+1}$ which fulfil the condition, differing in the value $\pi_{m+1}.$  For each case $\pi_2...\pi_{m+1}$ is represented by one of the permutations $\eta$ introduced above. However, the two cases $\pi$ with $|\pi_1-\pi_{m+1}|=1$ belong to the same $\eta .$ That is why we consider them together and divide their probability by two. Now
\[ p= \frac{P_m(\lambda)}{P_m(\kappa)}\cdot \sum_\eta P_m(\eta) = P_m(\lambda) \ .\]
This proves that $P_{m+1}$ extends $P_m.$ For the stationarity, the same proof has to be performed with extension to the left and condition (iii). We have to care that now $P_m(\kappa)$ refers to $\eta_2...\eta_m.$ However, since we assumed that $P_m$ is stationary, this is the same number, and the proof procedes as above.
\hfill$\Box$ \smallskip 

We called this a Markov extension since \eqref{ex} says that $\pi_{m+1}$ does not depend on $\pi_1.$ Since we do not assume any independence properties of $P_m,$ we cannot expect more. However, if we start with $P_2(12)=P_2(21)=\frac12$ and extend successively to $m=3,4,...$ we obtain the coin tossing order.

There are many other extensions. Just to give an example, we can divide the double cases in an unsymmetric way. A careful study of extensions may lead to a better model than coin tossing. But we have to stop here.

\section{Conclusion.} 
The six patterns of length 3 are naturally transformed into four independent contrasts. The turning rate $\alpha$ quantifies roughness and the parameters $\beta, \gamma, \delta$ describe deviations from symmetry and reversibility. For symmetric processes, to which EEGs seem to belong,  $\alpha$ controls the whole order structure of length 3. In an EEG, $\alpha$ is taken as function of time and delay, and its time oscillations seem to indicate an infra-slow rhythm of the brain. 

The next step is to study patterns of length 4. Theory and careful study of large datasets must go hand in hand. It was shown that any empirical distribution of patterns of length 4 extends to a distribution of patterns of arbitrary length. The main theoretical problem is the existence of self-similar distributions where pattern frequencies do not depend on delay.
Two contrasts $\tau_4$ and $\beta_4$ were suggested in a new paper by Weiss, Ruiz Marin, Keller and Matilla-Garcia \cite{MW} which continues the recent development of statistical tests with order patterns \cite{els,sou,wei}. The field is active and driven by new applications. It will form its own theory in the coming years.

\begin{acknowledgements} 
I am grateful to Karsten. And also to Jose Amigo and Osvaldo Rosso for the wonderful workshop in Dresden. I thank Bernd Pompe for discussions on my talk and this paper.
\end{acknowledgements}

\bigskip

\bigskip 

{\bf Supporting information:}  this work was not financially supported.

\end{document}